\documentclass[11pt,leqno]{article} \textwidth 160mm \textheight 235mm
\usepackage{amsfonts}
\usepackage{amssymb}
\usepackage{graphicx}

\def\ve{\varepsilon}

\def\emp{\emptyset}

\def\dom{{\rm dom}\,}

\def\N{{\cal N}}
\def\O{{\cal O}}

\def\R{I\!\!R}
\def\N{I\!\!N}

\def\ox{\overline{x}}
\def\oy{\overline{y}}
\def\oz{\overline{z}}

\def\disp{\displaystyle}

\def\Limsup{\mathop{{\rm Lim}\,{\rm sup}}}

\def\tto{\;{\lower 1pt \hbox{$\rightarrow$}}\kern -10pt
\hbox{\raise 2pt \hbox{$\rightarrow$}}\;}
\def\Hat{\widehat}
\def\Tilde{\widetilde}
\def\Bar{\overline}
\def\ra{\rangle}
\def\la{\langle}
\def\ve{\varepsilon}

\def\h{\hfill\Box}
\def\R{\mathbb{R}}
\def\N{I\!\!N}
\def\ox{\bar{x}}
\def\oy{\bar{y}}
\def\oz{\bar{z}}
\def\ov{\bar{v}}
\def\ow{\bar{w}}
\def\ou{\bar{u}}

\def\op{\bar{p}}

\def\ri{\mbox{\rm ri}\,}
\def\int{\mbox{\rm int}\,}
\def\gph{\mbox{\rm gph}\,}

\def\dim{\mbox{\rm dim}\,}
\def\dom{\mbox{\rm dom}\,}
\def\ker{\mbox{\rm ker}\,}

\def\aff{\mbox{\rm aff}\,}
\def\clco{\mbox{\rm clco}\,}

\def\substack#1#2{{\scriptstyle{#1}\atop\scriptstyle{#2}}}

\def\h{\hfill\triangle}
\def\dn{\downarrow}
\def\O{\Omega}
\def\ph{\varphi}
\def\emp{\emptyset}
\def\st{\stackrel}
\def\oR{\Bar{\R}}
\def\lm{\lambda}
\def\gg{\gamma}
\def\dd{\delta}
\def\al{\alpha}

\def \N{I\!\!N}
\def\th{\theta}
\def\sce{\setcounter{equation}{0}}
\begin{document}
\vspace*{0.7in}
\begin{center}
{\bf SECOND-ORDER SUBDIFFERENTIAL CALCULUS WITH\\ APPLICATIONS TO
TILT STABILITY IN OPTIMIZATION}
\\[2ex]
B. S. MORDUKHOVICH \footnote{Department of Mathematics, Wayne
State University, Detroit, MI 48202 (boris@math.wayne.edu).
Research of this author was partly supported by the National
Science Foundation under grant DMS-1007132, by
the Australian Research Council under grant DP-12092508, and by the Portuguese
Foundation of Science and Technologies under grant MAT/11109.} and R. T.
ROCKAFELLAR\footnote{Department of Mathematics, University of
Washington, Seattle, WA 98195 (rtr@math.washington.edu).}
\end{center}

\small{\bf Abstract.} The paper concerns the second-order
generalized differentiation theory of variational analysis and new
applications of this theory to some problems of constrained
optimization in finite-dimensional spaces. The main attention is
paid to the so-called (full and partial) second-order
subdifferentials of extended-real-valued functions, which are
dual-type constructions generated by coderivatives of first-order
subdifferential mappings. We develop an extended second-order
subdifferential calculus and analyze the basic second-order
qualification condition ensuring the fulfillment of the principal
second-order chain rule for strongly  and fully amenable compositions.
The calculus results obtained in this way and computing the second-order
subdifferentials for piecewise linear-quadratic functions and their major
specifications are applied then to the study of tilt stability
of local minimizers for important classes of problems in
constrained optimization that include, in particular, problems of
nonlinear programming and certain classes of extended nonlinear
programs described in composite terms.\vspace*{0.05in}

{\bf Key words.} variational analysis, constrained optimization,
nonlinear and extended nonlinear programming, second-order
subdifferentials, calculus rules, qualification conditions,
amenable functions, tilt-stable minimizers, strong regularity\vspace*{0.05in}

{\bf AMS subject classifications.} 49J52, 90C30, 90C31\vspace*{0.05in}

{\bf Abbreviated title.} Second-order subdifferential calculus

\newtheorem{Theorem}{Theorem}[section]
\newtheorem{Proposition}[Theorem]{Proposition}
\newtheorem{Remark}[Theorem]{Remark}
\newtheorem{Lemma}[Theorem]{Lemma}
\newtheorem{Corollary}[Theorem]{Corollary}
\newtheorem{Definition}[Theorem]{Definition}
\newtheorem{Example}[Theorem]{Example}
\renewcommand{\theequation}{{\thesection}.\arabic{equation}}
\renewcommand{\thefootnote}{\fnsymbol{footnote}}

\normalsize
\section{Introduction}\sce

Variational analysis has been recognized as a fruitful area of
mathematics, which primarily deals with optimization-related
problems while also applying variational principles and techniques
(largely based on perturbation and approximation ideas) to a broad
spectrum of problems that may not be of a variational nature. We
refer the reader to the books by Borwein and Zhu \cite{bz},
Mordukhovich \cite{m06a,m06b}, Rockafellar and Wets \cite{rw}, and
the bibliographies therein for the major results of variational
analysis and its numerous applications.

Since nonsmooth functions, sets with nonsmooth boundaries, and
set-valued mappings naturally and frequently appear in the
framework of variational theory and its applications via using
variational principles and techniques (even for problems with
smooth initial data), tools of generalized differentiation play a
crucial role in many aspects of variational analysis and
optimization; see, e.g., the books \cite{bz,c,dr-b,jl,m06a,m06b,rw,s} and
their references.

Over the years, the first-order subdifferential theory of
variational analysis has been well developed and understood in
both finite-dimensional and infinite-dimensional settings; see
\cite{bz,m06a,rw} and the commentaries therein. In
contrast, the {\em second-order} theory still requires a lot of
further development and implementation, although many
second-order generalized differential constructions have been
suggested and successfully applied to various optimization,
sensitivity, and related problems; see, e.g., the books
\cite{bs,m06a,rw} summarizing mainstream developments and trends
in the second-order theory and its applications.

As is well known, there are two generally independent approaches
to second-order differentiation in the classical analysis. One of
them is based on the Taylor expansion while the other defines the
second derivative of a function as the derivative of its
first-order derivative.\vspace*{0.05in}

In this paper we develop the latter ``derivative-of-derivative"
approach to the second-order generalized differentiation of
extended-real-valued functions
$\ph\colon\R^n\to\oR:=(-\infty,\infty]$ finite at the reference
points. The dual-space route in this vein suggested by
Mordukhovich \cite{m92} is to treat a (first-order) {\em
subdifferential} $\partial\ph$ of $\ph$ as a set-valued analog of
the classical derivative for nonsmooth functions and then to
define a {\em second-order subdifferential} $\partial^2\ph$ of
$\ph$ via a {\em coderivative} (generalized adjoint derivative
operator) $D^*\partial\ph$ of the subgradient mapping
$\partial\ph$; see Section~2 for more details. This second-order
construction was originally motivated by applications to
sensitivity analysis of variational systems \cite{m92,m94b}
inspired by the coderivative characterization of Lipschitzian
stability \cite{m92,m93}, but then the second-order
subdifferential and its modification were successfully employed in
the study of a broad spectrum of other important issues in
variational analysis and its applications; see, e.g.,
\cite{ccyy,dr,ew,hmn,hos,hr,jl,lm,lpr,m06a,m06b,mo,mo1,o,pr98,ye,z}
and the references therein. We specifically mention a remarkable result by Poliquin
and Rockafellar \cite{pr98} who established a full
characterization of {\em tilt-stable local minimizers} of functions (a new
notion introduced by them motivated, in particular, by the
justification of numerical algorithms) as the positive-definiteness of the second-order subdifferential mapping. For
${\cal C}^2$ functions, the latter criterion reduces to the
positive-definiteness of the classical Hessian matrix---a
well-known sufficient condition for the standard optimality in
unconstrained problems, which happens to be necessary and
sufficient for tilt-stable local minimizers \cite{pr98}. We also refer the reader
to the recent papers by Chieu et al. \cite{ccyy,ch} providing complete characterizations
of {\em convexity} and {\em strong convexity} of nonsmooth (in the second order)
functions via positive-semidefiniteness and definiteness of their second-order subdifferentials $\partial^2\ph$.
Related characterizations of {\em monotonicity} and {\em submonotonicity} of continuous mappings can be found in \cite{ct}.

Needless to say that efficient implementations and potential
extensions of the latter result to constrained optimization
problems, as well as any other valuable applications of the
aforementioned second-order subdifferential construction and its
modifications, largely depend on the possibility to develop a
fairly rich {\em second-order subdifferential calculus} and on
precisely {\em calculating} such constructions for attractive classes
of nonsmooth functions overwhelmingly encountered in variational
analysis and optimization. A certain amount of useful second-order
calculus rules were developed in \cite{lm,m94a,m02,m06a,mo,mo1,mw}. On the other hand,
precisely calculating the second-order subdifferential entirely in terms of the initial data
was effected for the following major classes of extended-real-valued functions particularly important in various applications:\vspace*{0.05in}

$\bullet$ For the indicator functions of {\em convex polyhedra} and related settings it was initiated by Dontchev and Rockafellar
\cite{dr} and then developed in \cite{bms,hmn,hos,hr,hy,nam,q,yy} for more involved frameworks. The obtained calculations played a crucial role in deriving \cite{dr}
verifiable characterizations of Robinson's strong regularity \cite{rob} for variational inequalities over (convex) polyhedral
sets as well as their specifications for complementarity problems and the associated Karush-Kuhn-Tucker (KKT) conditions for
nonlinear programs with ${\cal C}^2$ data. Further results in this vein on Lipschitzian stability of parametric variational
systems were given in \cite{bms,hmn,nam,q,yy1} and other publications in both finite and infinite dimensions. Applications to stationarity conditions for stochastic equilibrium
problems with equilibrium constraints in electricity spot market modeling were developed by Henrion and R\"omisch \cite{hr}.

$\bullet$ For the so-called {\em separable piecewise ${\cal C}^2$} functions it was done by Mordukhovich and Outrata \cite{mo}; see also \cite{ccyy} for further developments. it provided the basis for the efficient sensitivity analysis \cite{mo} of mathematical programs with equilibrium constraints (MPECs) including practical ones that arose in
applications to certain contact problems of continuum mechanics.

$\bullet$ For indicator functions to {\em smooth nonpolyhedral inequality systems} it was done by Henrion, Outrata, and Surowiec \cite{hos} by
employing and developing the transformation formulas from \cite{mo}. Then these calculations were applied in \cite{hos1,su} to deriving stationarity
conditions for equilibrium problems with equilibrium constraints (EPECs) in both deterministic and stochastic frameworks and to EPEC models of
oligopolistic competition in electricity spot markets.

$\bullet$ For the special class of functions arising in optimal control of the Moreau {\em sweeping process} is was done in the paper by Colombo et al. \cite{chhm}.
These calculations played a significant role in deriving constructive optimality conditions for discontinuous differential inclusions generated by the sweeping
process with great potentials for further applications. \vspace*{0.1in}

Now we briefly describe the main {\em goals and achievements} of this paper. Our primary attention is focused on the following
major issues new in second-order variational analysis:\vspace*{0.05in}

$\bullet$ Developing refined {\em second-order chain rules} of the
equality and inclusion (outer/upper estimate) types for
the aforementioned second-order subdifferential and its partial
modifications.

$\bullet$ Analyzing the basic {\em second-order qualification
condition} ensuring the fulfillment of the extended second-order
chain rules for {\em strongly amenable} compositions.

$\bullet$ Precise {\em calculating} second-order subgradients for major
classes of {\em fully amenable} functions.

$\bullet$ Applications of the obtained calculus and computational results to
deriving necessary optimality conditions as well as to
establishing {\em complete characterizations} of {\em tilt-stable
minimizers} for broad classes of constrained optimization problems
including those in nonlinear programming (NLP) and
extended nonlinear programming (ENLP) described via amenable compositions.\vspace{0.05in}

The rest of the paper is organized as follows. Section~2 contains
basic definitions  and brief discussions of the
first-order and second-order generalized differential
constructions studied and used in the paper. We also review there
some preliminary results widely employed in the sequel.

In Section~3 we deal with {\em second-order chain rules} of the
equality and inclusion types for the basic second-order subdifferential and its partial counterparts. The
equality-type results are established under the {\em full rank
condition} on the Jacobian matrix of the inner mapping of the
composition. Without imposing the latter assumption, we develop a
new {\em quadratic penalty} approach that allows us to derive
inclusion-type second-order chain rules for a broad class of
strongly amenable compositions valid under certain
second-order qualification conditions. The latter chain rules
generally provides merely outer estimates of the
second-order subgradient sets for compositions: we present an
example showing that the chain rule inclusion may be strict
even the linear inner mapping and piecewise linear outer functions
in fully amenable compositions.

Section~4 is devoted to a detailed analysis of the basic
{\em second-order qualification condition} ensuring the underlying
second-order chain rule for strongly amenable compositions.
Although the latter condition is automatically fulfilled under the
full rank assumption on the Jacobians of inner mappings as well as
for ${\cal C}^{1,1}$ outer functions in compositions, it seems to
be rather restrictive when outer functions are
extended-real-valued. In particular, we show that the second-order
qualification condition {\em reduces locally} to the {\em full
rank} requirement on the inner mapping Jacobian matrix if the
outer function is either convex piecewise linear, or it belongs to
a certain major class of piecewise linear-quadratic
functions. The results obtained in this direction are based on
precise {\em calculations} of the second-order subgradient sets of
the remarkable classes of fully amenable compositions under
consideration.

The concluding Section~5 concerns applications of the chain rules and
calculation results developed in the previous sections to the study of
{\em tilt-stable minimizers} for some classes of constrained optimization
problems represented in composite formats, which are convenient for developing
both theoretical and computational aspects of optimization. Such classes include,
besides standard nonlinear programs (NLP), broader models of the so-called extended
nonlinear programming (ENLP). Based on the second-order sum and chain rules with equalities,
we derive {\em complete characterizations} of tilt-stable local minimizers for important problems
of constrained optimization. The results obtained show, in particular, that for a general class of NLP problems the well-recognized
{\em strong second-order optimality condition} is {\em necessary and sufficient} for the tilt-stability
of local minimizers, which therefore is equivalent to Robinson's {\em strong regularity} of the associated variational
inequalities in such settings. Furthermore, the calculus rules obtained in this paper for {\em partial} second--order subdifferentials 
lead us also to characterizations of {\em full stability} in optimization (see Remark~\ref{fs}), while a detailed elaboration of this approach is a subject of our 
ongoing research. \vspace{0.05in}

Although a number of the results obtained in this paper hold in (or can be naturally extended to)
infinite-dimensional spaces, we confine ourselves for definiteness to the finite-dimensional setting.
Throughout the paper we use standard notation of variational analysis; cf.\ \cite{m06a,rw}. Recall that,
given a set-valued mapping $F\colon\R^n\tto\R^m$, the symbol
\begin{eqnarray}\label{1.1}
\begin{array}{ll}
\disp\Limsup_{x\to\ox}F(x):=\Big\{y\in\R^m\Big|&\exists\,
x_k\to\ox,\;\exists\,y_k\to y\;\mbox{ as }\;k\to\infty\\
&\mbox{with }\;y_k\in F(x_k)\;\mbox{ for all
}\;k\in\N:=\{1,2,\ldots\}\Big\}
\end{array}
\end{eqnarray}
signifies the {\em Painlev\'e-Kuratowski outer/upper limit} of $F$
as $x\to\ox$. Given a set $\O\subset\R^n$ and an
extended-real-valued function $\ph\colon\R^n\to\oR$ finite at
$\ox$, the symbols $x\st{\O}{\to}\ox$ and $x\st{\ph}{\to}\ox$
stand for $x\to\ox$ with $x\in\O$ and for $x\to\ox$ with
$\ph(x)\to\ph(\ox)$, respectively.

\section{Basic Definitions and Preliminaries}\sce

In this section we define and briefly discuss the basic
generalized differential constructions of our study and review
some preliminaries widely used in what follows; see \cite{m06a,rw} for more details.

Let $\ph\colon\R^n\to\oR$ be an extended-real-valued function
finite at $\ox$. The {\em regular subdifferential} (known also as
the presubdifferential and as the Fr\'echet or viscosity
subdifferential) of $\ph$ at $\ox$ is
\begin{eqnarray}\label{2.1}
\Hat\partial\ph(\ox):=\Big\{v\in\R^n\Big|\;\frac{\ph(x)-\ph(\ox)-\la
v,x-\ox\ra}{\|x-\ox\|}\ge 0\Big\}.
\end{eqnarray}
Each $v\in\Hat\partial\ph(\ox)$ is a {\em regular subgradient} of
$\ph$ at $\ox$. While $\Hat\partial\ph(\ox)$ reduces to a
singleton $\{\nabla\ph(\ox)\}$ if $\ph$ is Fr\'echet
differentiable at $\ox$ with the gradient $\nabla\ph(\ox)$ and to
the classical subdifferential of convex analysis if $\ph$ is
convex, the set (\ref{2.1}) may often be empty for nonconvex and
nonsmooth functions as, e.g., for $\ph(x)=-|x|$ at $\ox=0\in\R$.
Another serious disadvantage of the subdifferential construction
(\ref{2.1}) is the failure of standard calculus rules inevitably
required in the theory and applications of variational analysis
and optimization. In particular, the inclusion (outer estimate)
sum rule
$\Hat\partial(\ph_1+\ph_2)(\ox)\subset\Hat\partial\ph_1(\ox)+
\Hat\partial\ph_2(\ox)$
does not hold for the simplest nonsmooth functions $\ph_1(x)=|x|$
and $\ph_2(x)=-|x|$ at $\ox=0\in\R$.

The picture dramatically changes when we employ a limiting
``robust regularization" procedure over the subgradient mapping
$\Hat\partial\ph(\cdot)$ that leads us to the (basic first-order)
{\em subdifferential} of $\ph$ at $\ox$ defined by
\begin{eqnarray}\label{2.2}
\partial\ph(\ox):=\disp\Limsup_{x\st{\ph}{\to}\ox}\Hat\partial\ph(x)
\end{eqnarray}
and known also as the general, or limiting, or Mordukhovich
subdifferential; it was first introduced in \cite{m76} in an
equivalent way. Each $v\in\partial\ph(\ox)$ is called a (basic) {\em
subgradient} of $\ph$ at $\ox$. Thus, by taking into account
definition (\ref{1.1}) of $\Limsup$ and the notation
$x\st{\ph}{\to}\ox$, we represent the basic subgradients
$v\in\partial\ph(\ox)$ as follows:
\begin{eqnarray*}
\mbox{there are sequences }x_k\to\ox\;\mbox{ with
}\;\ph(x_k)\to\ph(\ox)\;\mbox{ and
}\;v_k\in\Hat\partial\ph(x_k)\;\mbox{ with }\;v_k\to v.
\end{eqnarray*}
In contrast to (\ref{2.1}), the subgradient set (\ref{2.2}) is
generally nonconvex (e.g., $\partial\ph(0)=\{-1,1\}$ for
$\ph(x)=-|x|$) while enjoying comprehensive calculus rules (``full
calculus"); this is based on {\em variational/extremal
principles}, which replace separation arguments in the absence of
convexity. Moreover, the basic subdifferential (\ref{2.2}) occurs
to be the {\em smallest} among any axiomatically defined
subgradient sets satisfying certain natural requirements; see
\cite[Theorem~9.7]{ms}.

In what follows we also need another subdifferential construction effective for non-Lipschitzian
extended-real-valued functions. Given $\ph\colon\R^n\to\oR$ finite at $\ox$, the {\em singular/horizontal
subdifferential} $\partial^\infty\ph(\ox)$ of $\ph$ at $\ox$ is defined by

\begin{equation}\label{sin}
\partial^\infty\ph(\ox):=\Limsup_\substack{x\st{\ph}{\to}\ox}{\lm\dn 0}\lm\Hat\partial\ph(x).
\end{equation}
If the function $\ph$ is lower semicontinuous (l.s.c.) around $\ox$, then $\partial^\infty\ph(\ox)=\{0\}$ if and only if $\ph$
is locally Lipschitzian around this point.

Given further a nonempty subset $\O\subset\R^n$, consider its
indicator function $\delta(x;\O)$ equal to 0 for
$x\in\O$ and to $\infty$ otherwise. For any fixed $\ox\in\O$,
define the {\em regular normal cone} to $\O$ at $\ox$ by
\begin{eqnarray}\label{2.3}
\Hat
N(\ox;\O):=\Hat\partial\delta(\ox;\O)=\disp\Big\{v\in\R^n\Big|\;\limsup_{x\st{\O}{\to}\ox}\frac{\la
v,x-\ox\ra}{\|x-\ox\|}\le 0\Big\}
\end{eqnarray}
and similarly the (basic, limiting) {\em normal cone} to $\O$ at
$\ox$ by $N(\ox;\O):=\partial\delta(\ox;\O)$. It follows from
(\ref{2.2}) and (\ref{2.3}) that the normal cone $N(\ox;\O)$
admits the limiting representation
\begin{eqnarray}\label{2.4}
N(\ox;\O)=\disp\Limsup_{x\st{\O}{\to}\ox}\Hat N(x;\O)
\end{eqnarray}
meaning that the basic normals $v\in N (\ox;\O)$ are those vectors
$v\in\R^n$ for which there are sequences $x_k\to\ox$ and $v_k\to
v$ with $x_k\in\O$ and $v_k\in\Hat N(x_k;\O)$, $k\in\N$. If $\O$
is locally closed around $\ox$, (\ref{2.4}) is equivalent to the
original definition by Mordukhovich \cite{m76}:
\begin{eqnarray*}
N(\ox;\O)=\disp\Limsup_{x\to\ox}\Big[\mbox{cone}\big(x-\Pi(x;\O)\big)\Big],
\end{eqnarray*}
where $\Pi(x;\O)$ signifies the Euclidean projector of
$x\in\R^n$ on the set $\O$, and where ``cone" stands for the conic
hull of a set.

There is the duality/polarity correspondence
\begin{eqnarray}\label{2.5}
\Hat N(\ox;\O)=T(\ox;\O)^*:=\Big\{v\in\R^n\Big|\;\la v,w\ra\le
0\;\mbox{ for all }\;w\in T(\ox;\O)\Big\}
\end{eqnarray}
between the regular normal cone (\ref{2.3}) and the {\em tangent
cone} to $\O$ at $\ox\in\O$ defined by
\begin{eqnarray}\label{2.6}
T(\ox;\O):=\Big\{w\in\R^n\Big|\;\exists\,x_k\in\O,\;\exists\,\al_k\ge
0\;\mbox{ with }\;\al_k(x_k-\ox)\to w\;\mbox{ as
}\;k\to\infty\Big\}
\end{eqnarray}
and known also as the Bouligand-Severi contingent cone to $\O$ at
this point. Note that the basic normal cone (\ref{2.4})
cannot be tangentially generated in a polar form (\ref{2.5}) by
using some set of tangents, since it is intrinsically
nonconvex while the polar $T^*$ to any set $T$ is
automatically convex. In what follows we may also use the subindex set notation
like $N_\O(\ox)$, $T_\O(\ox)$, etc. for the constructions involved.

It is worth observing that the {\em convex closure}
\begin{eqnarray}\label{2.7}
\Bar N(\ox;\O):=\clco N(\ox;\O)
\end{eqnarray}
of (\ref{2.4}), known as the Clarke/convexified normal cone to $\O$ at $\ox$ (see \cite{c}), may
dramatically enlarge the set of basic normals (\ref{2.4}). Indeed, it is proved by Rockafellar \cite{r85} that for every
vector function $f\colon\R^n\to\R^m$ locally Lipschitzian around
$\ox$ the convexified normal cone (\ref{2.7}) to the graph of $f$ at
$(\ox,f(\ox))$ is in fact a {\em linear subspace} of dimension
$d\ge m$ in $\R^n\times\R^m$, where the equality $d=m$ holds
if and only if the function $f$ is strictly differentiable
at $\ox$ with the derivative (Jacobian matrix) denoted for simplicity by $\nabla f(\ox)$, i.e.,
\begin{eqnarray*}
\disp\lim_{x,u\to\ox}\frac{f(x)-f(u)-\nabla
f(\ox)(x-u)}{\|x-u\|}=0,
\end{eqnarray*}
which is automatic when $f$ is ${\cal C}^1$ around $\ox$. In
particular, this implies that $\Bar N((\ox,f(\ox));\gph f)$ is the
whole space $\R^n\times\R^m$ whenever $f$ is nonsmooth
around $\ox$, $n=1$, and $m\ge 1$. Moreover, the aforementioned
results are discovered by Rockafellar \cite{r85} not only for
graphs of locally Lipschitzian functions, but also for the
so-called {\em Lipschitzian manifolds} (or graphically
Lipschitzian sets), which are locally homeomorphic to graphs of
Lipschitzian vector functions. The latter class includes graphs of
maximal monotone relations and subdifferential
mappings for convex, saddle, lower-${\cal C}^2$, and more general
prox-regular functions typically encountered in variational
analysis and optimization. In fact such {\em graphical sets} play a crucial role in the coderivative and
second-order subdifferential constructions studied in this paper.

Given a set-valued mapping $F\colon\R^n\tto\R^m$, define its {\em coderivative} at $(\ox,\oy)\in\gph F$ by \cite{m80}
\begin{equation}\label{cod}
D^* F(\ox,\oy)(v):=\Big\{u\in\R^n\big|\;(u,-v)\in N\big((\ox,\oy);\gph F\big)\Big\},\quad v\in\R^m,
\end{equation}
via the normal cone (\ref{2.4}) to the graph $\gph F$. Clearly the mapping $D^*F(\ox,\oy)\colon\R^m\tto\R^n$ is positive-homogeneous; it reduces to the adjoint derivative
\begin{equation}\label{smooth}
D^*F(\ox)(v)=\big\{\nabla F(\ox)^*v\big\},\quad v\in\R^m,
\end{equation}
where $^*$ stands for the matrix transposition, if $F$ is single-valued (then we omit $\oy=F(\ox)$ in the coderivative notation) and strictly differentiable at $\ox$. Note that the coderivative values in (\ref{cod}) are often nonconvex sets due to the nonconvexity of the normal cone on the right-hand side. Furthermore, the latter cone is taken to a graphical set, and thus its convexification in (\ref{cod}) may create serious troubles; see above. \vspace*{0.05in}

The main construction studied in the paper was introduced in \cite{m92} as follows.

\begin{Definition}{\bf (second-order subdifferential).}\label{2nd} Let the function $\ph\colon\R^n\to\oR$  be finite at $\ox$, and let $\oy\in\partial\ph(\ox)$  be a basic first-order subgradient of $\ph$ at $\ox$. Then the {\sc second-order subdifferential} of $\ph$ at $\ox$ relative to $\oy$ is defined by
\begin{equation}\label{2nd1}
\partial^2\ph(\ox,\oy)(u):=(D^*\partial\ph)(\ox,\oy)(u),\quad u\in\R^n,
\end{equation}
via the coderivative {\rm(\ref{cod})} of the first-order subdifferential mapping {\rm(\ref{2.2})}.
\end{Definition}

Observe that if $\ph\in{\cal C}^2$ around $\ox$ (in fact, it is merely continuous differentiable around $\ox$ with the strict differentiable first-order derivative at this point), then
\begin{eqnarray*}
\partial^2\ph(\ox)(u)=\big\{\nabla^2\ph(\ox)u\big\},\quad u\in\R^n,
\end{eqnarray*}
where $\nabla^2\ph(\ox)$ is the (symmetric) Hessian of $\ph$ at $\ox$. Sometimes the second-order construction (\ref{2nd1}) is called the ``generalized Hessian" of $\ph$ at the reference point \cite{pr98}. Note also that for the so-called {\em ${\cal C}^{1,1}$ functions} (i.e., continuously differentiable ones with locally Lipschitzian derivatives around $\ox$; another notation is ${\cal C}^{1+}$), we have the representation
\begin{eqnarray*}
\partial^2\ph(\ox)(u)=\partial\la u,\nabla\ph\ra(\ox),\quad u\in\R^n,
\end{eqnarray*}
via the basic first-order subdifferential (\ref{2.2}) of the derivative scalarization $\la u,\nabla\ph\ra(x):=\la u,\nabla\ph(x)\ra$ as $x\in\R^n$; see \cite[Proposition~1.120]{m06a}. It is worth emphasizing that the second-order subdifferential (\ref{2nd1}) as well as the generating coderivative and first-order subdifferential mappings
are {\em dual-space} intrinsically nonconvex constructions, which cannot correspond by duality to any derivative-like objects in primal spaces studied, e.g., in \cite{bs,rw}.

Following the scheme of Definition~\ref{2nd} and keeping the coderivative (\ref{cod}) as the underlying element of our approach while using different first-order subdifferentials in (\ref{2nd1}), we may define a variety of second-order constructions of type (\ref{2nd1}). In particular, for functions $\ph\colon\R^n\times\R^d\to\oR$ of $(x,w)\in\R^n\times\R^d$ there are two reasonable ways of introducing partial second-order subdifferentials; cf. \cite{lpr}. To proceed, define the {\em partial first-order} subgradient mapping $\partial_x\ph\colon\R^n\times\R^d\tto\R^n$ by
$$
\partial_x\ph(x,w):=\Big\{{\rm{set\;of\;subgradients}}\;v\;{\rm{of}}\;\ph_w:=\ph(\cdot,w)\;{\rm{at}}\;x\Big\}=\partial\ph_w(x).
$$
Then given $(\ox,\ow)$ and $\oy\in\partial_x\ph(\ox,\ow)$, define the {\em partial second-order subdifferential} of $\ph$ with respect to $x$ of at $(\ox,\ow)$ relative to $\oy$ by
\begin{equation}\label{par1}
\partial^2_x\ph(\ox,\ow,\oy)(u):=(D^*\partial\ph_{\ow})(\ox,\oy)(u)=\partial^2\ph_{\ow}(\ox,\oy)(u),\quad u\in\R^n,
\end{equation}
with $\ph_{\ow}(x)=f(x,\ow)$. On the other hand, we can define the {\em extended  partial second-order subdifferential} of $\ph$ with respect to $x$ of at $(\ox,\ow)$ relative to $\oy$ by
\begin{equation}\label{par2}
\Tilde\partial^2_x\ph(\ox,\ow,\oy)(u):=(D^*\partial_x\ph)(\ox,\ow,\oy)(u),\quad u\in\R^n.
\end{equation}
As argued in \cite{lpr}, constructions (\ref{par1}) and (\ref{par2}) are not the same even in the case of ${\cal C}^2$ functions  when (\ref{par1}) reduces to $\nabla^2_{xx}\ph(\ox,\ow)(u)$ while (\ref{par2}) comes out as $(\nabla^2_{xx}\ph(\ox,\ow)u,\nabla^2_{xw}\ph(\ox,\ow)u)$. This happens due to the involvement of $w\to\ow$ in the limiting procedure to define the extended partial second-order subdifferential set $\Tilde\partial^2_x\ph(\ox,\ow,\oy)(u)$, which is  hence larger than (\ref{par1}). Note that both partial second-order constructions (\ref{par1}) and (\ref{par2}) are proved to be useful in applications; see, e.g., \cite{lm,lpr} for more details.

It has been well recognized and documented (see, e.g., \cite{bz,m06a,m06b,rw,s} and the references therein) that the first-order limiting constructions (\ref{2.2}), (\ref{2.4}), and (\ref{cod}) enjoy full calculi, which are crucial for their numerous applications. Based on definitions (\ref{2nd1}) of the second-order subdifferential and its partial counterparts (\ref{par1}) and (\ref{par2}), it is natural to try to combine calculus results for first-order subgradients with those for coderivatives to arrive at the corresponding second-order calculus rules. However, there are nontrivial complications to proceed in this way due to the fact that general results of the first-order subdifferential calculus hold as {\em inclusions} while the coderivative (\ref{cod}) does {\em not} possesses any {\em monotonicity} properties. Thus the initial requirement arises on selecting classes of functions for which calculus rules for first-order subgradients hold as {\em equalities}. Proceeding in this direction, a number of second-order calculus rules
have been established in \cite{lm,m94a,m02,m06a,mo,mo1,mw} for full while {\em not for partial} second-order subdifferentials in finite and infinite dimensions.\vspace*{0.05in}

In the next section we obtain new second-order chain rules applied to full and partial second-order subdifferentials and develop, in particular, a direct approach based on quadratic penalties to derive general results for strongly amenable compositions.

\section{Second-Order Subdifferential Chain Rules}\sce

Given a vector function $h\colon\R^n\times\R^d\to\R^m$ with $m\le n$ and a proper extended-real-valued function $\th\colon\R^m\to\oR$, consider the composition
\begin{equation}\label{comp}
\ph(x,w)=(\th\circ h)(x,w):=\th\big(h(x,w)\big)
\end{equation}
with $x\in\R^n$ and $w\in\R^d$. Our first theorem provides {\em exact formulas} for calculating the partial second-order subdifferentials (\ref{par1}) and (\ref{par2}) of composition (\ref{comp}) under the {\em full rank condition} on the partial derivative (Jacobian matrix) $\nabla_x h(\ox,\ow)$ at the reference point.

\begin{Theorem}{\bf (exact second-order chain rules with full rank condition).}\label{full} Given a point $(\ox,\ow)\in\R^n\times\R^d$, suppose that $\th$ in {\rm(\ref{comp})} is finite at $\oz:=h(\ox,\ow)$, that $h(\cdot,\ow)\colon\R^n\to\R^m$ is continuously differentiable around $\ox$ with the full row rank condition
\begin{equation}\label{rank}
{\rm rank}\,\nabla_x h(\ox,\ow)=m,
\end{equation}
and that the mapping $\nabla_x h(\cdot,\ow)\colon\R^n\to\R^m$ is strictly differentiable at $\ox$. Pick any $\oy\in\partial_x\ph(\ox,\ow)$ and denote by $\ov$ a unique vector satisfying the relationships
$$
\ov\in\partial\th(\oz)\;\mbox{ and }\;\nabla_x h(\ox,\ow)^*\ov=\oy.
$$
Then we have the chain rule equality for the partial second-order subdifferential {\rm(\ref{par1})}:
\begin{equation}\label{rank1}
\partial^2_x\ph(\ox,\ow,\oy)(u)=\nabla^2_{xx}\la\ov,h\ra(\ox,\ow)u+\nabla_x h(\ox,\ow)^*\partial^2\th(\oz,\ov)(\nabla_x h(\ox,\ow)u),\quad u\in\R^n.
\end{equation}
If in addition the mapping $h\colon\R^n\times\R^d\to\R^m$ is continuously differentiable around $(\ox,\ow)$ and its derivative $\nabla h\colon\R^n\times\R^d\to\R^m$ is strictly differentiable at $(\ox,\ow)$, then we have
\begin{equation}\label{rank2}
\begin{array}{ll}
\Tilde\partial^2_x\ph(\ox,\ow,\oy)(u)&=\Big(\nabla^2_{xx}\la\ov,h\ra(\ox,\ow)u,\nabla^2_{xw}\la\ov,h\ra(\ox,\ow)u\Big)\\\\
&+\Big(\nabla_x h(\ox,\ow)^*\partial^2\th(\oz,\ov)(\nabla_x h(\ox,\ow)u),\nabla_w h(\ox,\ow)^*\partial^2\th(\oz,\ov)(\nabla_x h(\ox,\ow)u)\Big)
\end{array}
\end{equation}
whenever $u\in\R^n$ for the extended partial second-order subdifferential {\rm(\ref{par2})}.
\end{Theorem}
{\bf Proof.} We derive the chain rule (\ref{rank2}) for the extended partial second-order subdifferential; the proof of (\ref{rank1}) is just a simplification of the one given below.

On the first-order subdifferential level we have from \cite[Proposition~1.112]{m06a} under the assumptions made (and from \cite[Exercise~10.7]{rw} under some additional assumptions) that there is a neighborhood $U$ of $(\ox,\ow)$ such that
\begin{eqnarray*}
\partial_x\ph(x,w)=\Big\{y\in\R^n\Big|\;\exists\,v\in\partial\th\big(h(x,w)\big)\;\mbox{ with }\;\nabla_x h(x,w)^*v=y\Big\}\;\mbox{ for all }\;(x,w)\in U.
\end{eqnarray*}
For any fixed $\oy\in\partial_x\ph(\ox,\ow)$, this gives us locally around $(\ox,\ow,\oy)$ the graph representation
\begin{equation}\label{gr1}
\begin{array}{ll}
\gph\partial_x\ph=\Big\{(x,w,y)\in\R^n\times\R^d\times\R^n\Big|&\exists\,(p,v)\in\gph\th\;\mbox{ such that}\\
&h(x,w)=p,\;\nabla_x h(x,w)^*v=y\Big\}.
\end{array}
\end{equation}
Consider now the two possible cases in the graph representation (\ref{gr1}): {\bf (i)} the ``square" case when $m=n$ and {\bf (ii)} the ``general" one when $m<n$.

In the square case (i) we have by the full rank condition (\ref{rank}) that the matrix $\nabla_x h(x,w)$ is invertible for $(x,w)$ near $(\ox,\ow)$, and hence (\ref{gr1}) can be rewritten as
\begin{equation}\label{gr2}
\gph\partial_x\ph=\Big\{(x,w,y)\in\R^n\times\R^d\times\R^n\Big|\;f(x,w,y)\in\gph\partial\th\Big\}
\end{equation}
via the mapping $f\colon\R^n\times\R^d\times\R^n\to\R^{2n}$ given by
\begin{equation}\label{f}
f(x,w,y):=\Big(h(x,w),(\nabla_x h(x,w)^*)^{-1}y\Big)\;\mbox{ for }\;(x,w,y)\;\mbox{ near }\;(\ox,\ow,\oy).
\end{equation}
In other words, representation (\ref{gr2}) can be expressed via the preimage/inverse image of the set $\gph\partial\th$ under the mapping $f$ as follows:
\begin{equation}\label{pre}
\gph\partial_x\ph=f^{-1}\big(\gph\partial\th\big).
\end{equation}
Since $\nabla_x h$ is assumed to be strictly differentiable at $(\ox,\ow)$, the mapping $f$ in (\ref{f}) is strictly differentiable at $(\ox,\ow,\oy)$ and, by (\ref{rank}) with $m=n$, its Jacobian matrix $\nabla f(\ox,\ow,\oy)$ has full row rank $2n$. Employing \cite[Theorem~1.17]{m06a} to (\ref{pre}) gives us
\begin{equation}\label{sur}
N\big((\ox,\ow,\oy);\gph\partial_x\ph\big)=\nabla f(\ox,\ow,\oy)^*N\big(f(\ox,\ow,\oy);\gph\partial\th\big).
\end{equation}

Now we calculate the derivative/Jacobian matrix of $f$ at $(\ox,\ow,\oy)$ by using the particular structure of $f$ in (\ref{f}), the classical chain rule, and the well-known Leach inverse function theorem for strictly differentiable mappings; see, e.g., \cite{dr-b,m06a}. Define the mappings $f_1\colon\R^n\times\R^d\times\R^n\to\R^n$ and $f_2\colon\R^n\times\R^d\times\R^n\to\R^n$ with $f=(f_1,f_2)$ by $f_1(x,w,y):=f(x,w)$ and
$$
f_2(x,w,y):=(\nabla_x h(x,w)^{-1})^*y\;\mbox{ for }\;(x,w,y)\in\R^n\times\R^d\times\R^n.
$$
It is clear that $\nabla f_1(\ox,\ow,\oy)=(\nabla f(\ox,\bar w),0)$, while for calculating $\nabla f_2(\ox,\ow,\oy)$ we introduce two auxiliary mapping $g\colon\R^n\times\R^d\times\R^n\to\R^n$ and $q\colon\R^n\times\R^d\times\R^n\to\R$ by
\begin{eqnarray*}
g(x,w,p):=\nabla_x h(x,w)^*p\;\mbox{ and }\;q(x,w,p):=\la p,h(x,w)\ra\;\mbox{ for }\;(x,w,p)\in\R^n\times\R^d\times\R^n.
\end{eqnarray*}
Note that $g(x,w,p)=\nabla q_x (x,w,p)^*$ and that $g(x,w,f_2(x,w,y))-v=0$. Differentiating the latter equality gives us
\begin{equation}\label{g}
\nabla_x g(x,w,f_2(x,w,y))+\nabla_p g(x,w,f_2(x,w,y))\nabla_x f_2(x,w,y)=0.
\end{equation}
Observing further that $\nabla_x g(x,w,p)=\nabla_x(\nabla_x q(x,w,p))^*=\nabla_{xx}^2 q(x,w,p)$, we get from (\ref{g}) and the definitions above that the equation
\begin{eqnarray*}
\nabla^2_{xx}\la(\nabla_x h(x,w)^{-1})^*y,h(x,w)\ra+\nabla_x h(x,w)^*\nabla f_2(x,w,p)=0
\end{eqnarray*}
is satisfied, which implies in turn the representation of the partial derivative
\begin{equation}\label{x}
\nabla_x f_2(x,w,y)=-(\nabla_x h(x,w)^{-1})^*\cdot\nabla^2_{xx}\la(\nabla_x h(x,w)^{-1})^*y,h(x,w)\ra.
\end{equation}
Similarly we have the following representation of the other partial derivative of $f_2$:
\begin{equation}\label{w}
\nabla_w f_2(x,w,y)=-(\nabla_x h(x,w)^{-1})^*\cdot\nabla^2_{xw}\la(\nabla_x h(x,w)^{-1})^*y,h(x,w)\ra.
\end{equation}
Taking into account that $\nabla_y f_2(x,w,y)=(\nabla_x h(x,w))^{-1})^*$ gives us finally
\begin{eqnarray}\label{lea}
\nabla f(\ox,\bar w,\bar y)=\left[\begin{array}{clcr}
\nabla_x h(\ox,\bar w)&\nabla_w h(\ox,\bar w)&0\\
\quad\nabla_x f_2(\ox,\ow,\oy)&\nabla_w f_2(\ox,\ow,\oy)&(\nabla_x h(\ox,\bar w)^{-1})^*
\end{array}\right],
\end{eqnarray}
where $\nabla_x f_2(\ox,\ow,\oy)$ and $\nabla_w f_2(\ox,\ow,\oy)$ are calculated in (\ref{x}) and (\ref{w}), respectively. The second-order chain rule (\ref{rank2}) in the square case (i) follows now from substituting (\ref{lea}) into (\ref{sur}) and by using then the definitions of the constructions involved and elementary transformations.\vspace*{0.05in}

It remains to consider the general case (ii) with $m<n$. This case can be reduced to the previous one by introducing a linear mapping $\Tilde h\colon\R^n\times\R^d\to\R^{n-m}$ such that the mapping
$$
\Bar h(x,w):=\big(h(x,w),\Tilde h(x,w))\;\mbox{ from }\;\R^n\times\R^d\;\mbox{ to }\;\R^n
$$
has full rank. It can be done, e.g., by choosing a basis $\{a_1,\ldots,a_{n-m}\}$ for the ${n-m}$-dimensional spaces $\{u\in\R^n|\;\nabla h_x(\ox,\ow)u=0\}$ and letting
$\Bar h(x,w):=(h(x,w),\la a_1,x\ra,\ldots,\la a_{n-m},x\ra)$; cf.\ \cite[Exercise~6.7]{rw} for a first-order setting. Then viewing $\ph$ as $\Bar\th\circ\Bar h$ with $\Bar\th(z,p):=\th(z)$ for all $z\in\R^m$ and $p\in\R^{n-m}$ reduces (ii) to (i) and thus completes the proof of the theorem. $\h$\vspace*{0.05in}

Some remarks on the results related to those obtained in Theorem~\ref{full} are in order.

\begin{Remark}{\bf (discussions on second-order chain rules with full rank/surjectivity conditions).} {\rm Previously known second-order chain rules of type (\ref{rank1}) were derived for the full second-order subdifferential (\ref{2nd1}), where condition (\ref{rank}) was written as ${\rm rank}\nabla h(\ox)=m$. To the best of our knowledge, the first result in this direction was obtained in \cite[Theorem~3.4]{mo} with the inclusion ``$\subset$" in (\ref{rank1}). Various infinite-dimensional extensions of (\ref{rank1}) in the inclusion and equality forms were derived in \cite{m02,mw} and \cite[Theorem~1.127]{m06a} with imposing the {\em surjectivity condition} on the derivative $\nabla h(\ox)$ as the counterpart of (\ref{rank}) in infinite-dimensional spaces. Observe that the proof of (\ref{rank2}) given above in case (i) corresponding to the invertible partial derivative $\nabla_x h(\ox,\ow)$ holds in any Banach space, while the device in case (ii) is finite-dimensional.}
\end{Remark}

Next we explore the possibility of deriving second-order chain rules for (\ref{comp}) when the rank condition (\ref{rank}) may not be satisfied. This can be done for broad classes of amenable functions defined in the way originated in \cite{pr92}, which are overwhelmingly encountered in finite-dimensional parametric optimization. Recall \cite{lr} that a proper function $\ph\colon\R^n\times\R^d\to\oR$ is {\em strongly amenable} in $x$ at $\ox$ with {\em compatible parametrization} in $w$ at $\ow$ if there is a neighborhood $V$ of $(\ox,\ow)$ on which $\ph$ is represented in the composition form (\ref{comp}), where $h$ is of class ${\cal C}^2$ while $\th$ is a proper, l.s.c., convex function such that the first-order qualification condition
\begin{equation}\label{1qc}
\partial^\infty\th\big(h(\ox,\ow)\big)\cap\ker\nabla_x h(\ox,\ow)^*=\{0\}
\end{equation}
involving the singular subdifferential (\ref{sin}) is satisfied. The latter qualification condition automatically holds if either $\th$ is locally Lipschitzian around $h(\ox,\ow)$ or the full rank condition (\ref{rank}) is fulfilled, since it is equivalent to
\begin{eqnarray*}
\ker\nabla_x h(\ox,\ow)^*:=\big\{v\in\R^n\big|\;0=\nabla_x h(\ox,\ow)^*v\}=\{0\}.
\end{eqnarray*}
Properties of strongly amenable compositions $\ph(x)=\th(h(x))$ and related functions are largely investigated in \cite{pr92,pr96,rw}; most of them hold also for strongly amenable compositions (\ref{comp}) with compatible parametrization \cite{lpr,lr}. Strong amenability is a property that bridges between smoothness and convexity while covering at the same time a great many of the functions that are of interest as the essential objective in minimization problems; see \cite{rw} for more details.\vspace*{0.05in}

The next theorem establishes second-order subdifferential chain rules of the inclusion type for strongly amenable compositions with no full rank requirement (\ref{rank}).

\begin{Theorem}{\bf (second-order chain rules for strongly amenable compositions).}\label{am} Let $\ph\colon\R^n\times\R^d\to\oR$ be strongly amenable in $x$ at $\ox$ with compatible parametrization in $w$ at $\ow$, and let $\oy\in\partial_x\ph(\ox,\ow)$. Denote $\oz:=h(\ox,\ow)$ and consider the nonempty set
$$
M(\ox,\ow,\oy):=\Big\{v\in\R^m\Big|\;v\in\partial\th(\oz)\;\mbox{ with }\;\nabla_x h(\ox,\ow)^*v=\oy\Big\}.
$$
Assume the fulfillment of the second-order qualification condition:
\begin{eqnarray}\label{2qc}
\partial^2\th(\oz;v)(0)\cap\ker\nabla_x h(\ox,\ow)^*=\{0\}\;\mbox{ for all }\;v\in M(\ox,\ow,\oy).
\end{eqnarray}
Then we have the following chain rules for the partial second-order subdifferentials {\rm(\ref{par1})} and {\rm(\ref{par2})}, respectively, valued for all $u\in\R^n$:
\begin{equation}\label{am1}
\partial^2_x\ph(\ox,\ow,\oy)(u)\subset\bigcup_{v\in M(\ox,\ow,\oy)}\nabla^2_{xx}\la v,h\ra(\ox,\ow)u+\nabla_x h(\ox,\ow)^*\partial^2\th(\oz, v)(\nabla_x h(\ox,\ow)u),
\end{equation}
\begin{equation}\label{am2}
\begin{array}{ll}
\Tilde\partial^2_x\ph(\ox,\ow,\oy)(u)&\subset\disp\bigcup_{v\in M(\ox,\ow,\oy)}\Big(\nabla^2_{xx}\la v,h\ra(\ox,\ow)u,\nabla^2_{xw}\la v,h\ra(\ox,\ow)u\Big)\\\\
&+\Big(\nabla_x h(\ox,\ow)^*\partial^2\th(\oz, v)(\nabla_x h(\ox,\ow)u),\nabla_w h(\ox,\ow)^*\partial^2\th(\oz,v)(\nabla_x h(\ox,\ow)u)\Big).
\end{array}
\end{equation}
\end{Theorem}
{\bf Proof.} For brevity and simplicity of the arguments and notation, we present a detailed proof just for the full second-order subdifferential (\ref{2nd1}) of the strongly amenable nonparameterized compositions $\ph(x)=\th(h(x))$ in which case both formulas (\ref{am1}) and (\ref{am2}) reduce to
\begin{equation}\label{am3}
\partial^2\ph(\ox,\oy)(u)\disp\subset\bigcup_\substack{v\in\partial\th(\oz)}{\nabla h(\ox)^*v=\oy}\Big(\nabla^2\la v,h\ra(\ox)u+\nabla h(\ox)^*\partial^2\th(\oz,v)(\nabla h(\ox)u)\Big)
\end{equation}
with $\oz=h(\ox)$ under the  {\em basic second-order qualification condition}
\begin{equation}\label{2qc1}
\partial^2\th(\oz;v)(0)\cap\ker\nabla h(\ox)^*=\{0\}\;\mbox{ whenever }\;v\in\partial\th(\oz)\;\mbox{ and }\;\nabla h(\ox)^*v=\oy.
\end{equation}
The reader can readily check that the method of quadratic penalties developed below perfectly works for the case of partial second-order subdifferentials to produce the chain rule inclusions (\ref{am1}) and (\ref{am2}) under the ``partial" second-order qualification condition (\ref{2qc}).

We begin with observing that the first-order chain rule
\begin{equation}\label{1-am}
\partial\ph(x)=\nabla h(x)^*\partial\th(h(x))\;\mbox{ whenever }\;x\in U
\end{equation}
holds as equality for strongly amenable compositions on some neighborhood $U$ of $\ox$. Indeed, it follows from the more general chain rule of \cite[Theorem~3.41(iii)]{m06a} due to (\ref{smooth}) and the particular properties of strongly amenable functions summarized in \cite[Exercise~10.25]{rw}.

Now we proceed with calculating of the second-order subdifferential $\partial^2\ph(\ox,\oy)$ for the given first-order subgradient $\oy\in\partial\ph(\ox)$. The definitions in (\ref{2nd}), (\ref{cod}), and (\ref{2.4}) suggest to us calculating the regular normal cone $\Hat N((x,y);\O)$ to the subdifferential graph $\O:=\gph\partial\ph$ of $\ph$ at points $(x,y)\in\gph\partial\ph$ near $(\ox,\oy)$ and then passing to the limit therein as $(x,y)\to(\ox,\oy)$. To simplify notation, let us focus first on calculating $\Hat N((\ox,\oy);\O)$ for the graphical set $\O$. Developing a {\em variational approach} to subdifferential calculus and employing the {\em smooth variational description} of regular normals from \cite[Theorem~6.11]{rw} and \cite[Theorem~1.30]{m06a}, we have that $(\omega,-\xi)\in\Hat N((\ox,\oy);\O)$ if and only if there is a smooth function $\vartheta\colon\R^n\times\R^n\to\R$ such that
\begin{equation}\label{svd}
\disp{\rm argmin}_{(x,y)\in\O}\,\vartheta(x,y)=\{(\ox,\oy)\}\;\mbox{ and }\;\nabla\vartheta(\ox,\oy)=(-\omega,\xi).
\end{equation}
Using the first-order chain rule formula (\ref{1-am}) allows us to transform the minimization problem in (\ref{svd}) into the following one:
\begin{eqnarray}\label{min}
\left\{\begin{array}{ll}
\mbox{minimize }\vartheta\Big(x,\nabla h(x)^*v\Big)\;\mbox{ over all}\\
x\in\R^n,\;(v,z)\in\gph\partial\th\;\mbox{ with }\;h(x)-z=0.
\end{array}\right.
\end{eqnarray}
We know from (\ref{svd}) that $(x,z,v)$ is an optimal solution to (\ref{min}) if and only if
\begin{equation}\label{min1}
x=\ox,\;z=\oz,\;\mbox{ and }\;\nabla h(\ox)^*v=\oy.
\end{equation}
Let $G:=\gph\th$ and for any $\ve>0$ consider the {\em quadratic penalty problem}:
\begin{eqnarray}\label{min2}
\left\{\begin{array}{ll}
\mbox{minimize }\vartheta\Big(x,\nabla h(x)^*v\Big)+\disp\frac{1}{2\ve}\|h(x)-z\|^2\\
\mbox{over all }\;x\in\R^n\;\mbox{ and }\;(z,v)\in G.
\end{array}\right.
\end{eqnarray}
Denoting by $\ve\mapsto S(\ve)$ the optimal solution map for problem (\ref{min2}), observe that it is of {\em closed graph} and uniformly {\em bounded} around $\ve=0$. Indeed. it follows from the closedness of $\gph\partial\th$ due to the convexity and lower semicontinuity of $\th$, the local continuity of $\vartheta$ and $\nabla h$, and the uniform boundedness near $(\ox,\oy)$ of the mapping $M\colon\O\tto\R^m$ given by
\begin{eqnarray}\label{M}
M(x,y):=\Big\{v\in\partial\th(h(x))\Big|\;\nabla h(x)^*v=y\Big\}
\end{eqnarray}
with $\O=\gph\partial\ph$. The latter local boundedness can be easily verified arguing by contradiction due the qualification condition (\ref{1qc}) in the definition of amenable functions.

Now consider a sequence of values $\ve_k\to 0$ as $k\to\infty$ and pick a triple $(x_k,z_k,v_k)\in S(\ve_k)$ for all $k\in\N$. By the local boundedness of the solution map $S(\ve)$ near $\ve=0$ the sequence $\{v_k\}$ is bounded, and thus it has a cluster point $\ov$. Without loss of generality we suppose that $v_k\to\ov$ as $k\to\infty$ and get therefore that
\begin{equation}\label{lim}
x_k\to\ox,\;z_k:=h(x_k)\to\oz,\;\mbox{ and }\;y_k:=\nabla h(x_k)^*v_k\to\nabla h(\ox)^*\ov=\oy
\end{equation}
with $(\oz,\ov)\in G$. By the observation above (\ref{min1}), the triple $(\ox,\oz,\ov)$ is the unique optimal solution to the unperturbed problem (\ref{min}). On the other hand, applying the first-order necessary optimality conditions from \cite[Theorem~6.12]{rw} to the solution $(x_k,z_k,v_k)$ of the penalized problem (\ref{min2}) with a smooth cost function and a geometric constraint gives us
$$
\nabla_x\Big[\vartheta\Big(x,\nabla h(z)^*v_k\Big)+\frac{1}{\ve_k}\|h(x)-z_k\|^2\Big]\Big|_{x=x_k}=0,
$$
$$
-\nabla_{z,v}\Big[\vartheta\Big(x_k,\nabla h(x_k)^*v\Big)+\frac{1}{2\ve_k}\|h(x_k)-z\|^2\Big]\Big|_{(z,v)=(z_k,v_k)}\in\Hat N((z_k,v_k);G)
$$
for all $k\in\N$. Denoting $p_k:=[h(x_k)-z_k]/\ve_k$, these conditions calculate out to
\begin{equation}\label{op1}
\nabla_x\vartheta(x_k,y_k)+\nabla^2\la v_k,h\ra(x_k)\nabla_y\vartheta(x_k,y_k)+\nabla\la p_k,h\ra(x_k)=0,
\end{equation}
\begin{equation}\label{op2}
\Big(p_k,-\nabla h(x_k)\nabla_y\vartheta(x_k,y_k)\Big)\in\Hat N((z_k,v_k);G).
\end{equation}
By passing above to subsequences as $k\to\infty$ if needed, we can reduce the situation to considering one of the following two cases:\\[1ex]
{\bf Case~1:} $\{p_k\}$ converges to some $\op$.\\
{\bf Case~2:} $p_k\to\infty$ while $\{p_k/\|p_k\|\}$ converges to some $\op\ne 0$. \vspace*{0.1in}

In Case~1 it follows from (\ref{lim}), (\ref{op1}), and (\ref{op2}) that
\begin{equation}\label{op3}
\nabla_x\vartheta(\ox,\oy)+\nabla^2\la\ov,h\ra(\ox)\nabla_y\vartheta(\ox,\oy)+\nabla\la\op,h\ra(\ox)=0,
\end{equation}
\begin{equation}\label{op4}
\Big(\op,-\nabla h(\ox)\nabla_y\vartheta(\ox,\oy)\Big)\in\Hat N((\oz,\ov);G),
\end{equation}
where $\nabla_x\vartheta(\ox,\oy)=-\omega$ and $\nabla_y\vartheta(\ox,\oy)=\xi$ by the second equality in (\ref{svd}).

In Case~2 we get, dividing first both parts of (\ref{op1}) and (\ref{op2}) by $\|p_k\|$ and then passing to the limit therein as $k\to\infty$, that
\begin{equation}\label{op5}
\nabla h(\ox)^*\op=0\;\mbox{ and }\;(\op,0)\in\Hat N((\oz,\ov);G)\;\mbox{ with }\;\|\op\|=1.
\end{equation}
Thus, by talking into account our choice of $(\omega,-\xi)\in\Hat N((\ox,\oy);\O)$ and the construction of $M$ in (\ref{M}), we deduce from (\ref{op3})--(\ref{op5}) the existence of $\ov\in M(\ox,\oy)$ and $\op$ satisfying either (\ref{op3}) and (\ref{op4}) or (\ref{op5}). Since the arguments above equally hold for every point $(x,y)\in\gph\partial\ph$ near $(\ox,\oy)$, they ensure the following description of the regular normal cone to $\O=\gph\ph$ at points $(x,y)\in\O$ in a neighborhood of the reference one $(\ox,\oy)$, where $G=\gph\partial\th$:
\begin{equation}\label{desc}
\left[\begin{array}{ll}
(\omega,-\xi)\in\Hat N((x,y);\O)\Longrightarrow\;\exists\,v\in M(x,y),\;p\in\R^m\;\mbox{ such that}\\\\
\mbox{either: }\;\left\{\begin{array}{ll}
\omega=\nabla^2\la v,h\ra(x)\xi+\nabla h(x)^*p\\
\mbox{with }\;(p,-\nabla h(x)\xi)\in\Hat N((h(x),v);G)
\end{array}\right.\\\\
\mbox{or: }\;\nabla h(x)^*p=0\;\mbox{ with }\;(p,0)\in\Hat N(h(x),v);G),\;\|p\|=1.
\end{array}\right.
\end{equation}

Next we take any basic normal $(\omega,-\xi)\in N(\ox,\oy);\O)$ (not just a regular one) and by (\ref{2.4}) find sequences $(x_k,y_k)\to(\ox,\oy)$ and $(\omega_k,-\xi_k)\to(\omega,-\xi)$ as $k\to\infty$ satisfying
$$
(x_k,y_k)\in\O\;\mbox{ and }\;(\omega_k,-\xi_k)\in\Hat N((x_k,y_k);\O)\;\mbox{ for all }\;k\in\N.
$$
Employing the description of regular normals (\ref{desc}) ensures the existence of $v_k\in M(x_k,y_k)$ and $p_k\in\R^m$ such that the either/or alternative in (\ref{desc}) holds for each $k\in\N$. Due to the established local boundedness of the mapping $M$, suppose with no loss of generality that
$$
v_k\to v\;\mbox{ as }\;k\to\infty\;\mbox{ for some }\;v\in M(\ox,\oy).
$$
By a further passage to subsequences, we can reduce the situation to where just one of the "either/or" parts of the alternative in (\ref{desc}) holds for all $k$. Consider first the ``or" part of this alternative, i.e., the validity of
$$
\nabla h(x_k)^*p_k=0\;\mbox{ with }\;(p_k,0)\in\Hat N(h(x_k),v_k);G),\;\|p_k\|=1\;\mbox{ for all }\;k.
$$
In this case the sequence $\{p_k\}$ has a cluster point $p$, and thus we get
\begin{equation}\label{desc1}
\nabla h(\ox)^*p=0\;\mbox{ with }\;(p,0)\in N((h(\ox),v);G),\;\|p\|=1,\;\mbox{ and }\;v\in M(\ox,\oy)
\end{equation}
by passing to the limit as $k\to\infty$ and taking into account (\ref{2.4}) and the continuity of $f$ and $\nabla h$.

When the ``either" part holds, we proceed similarly to Cases~1 and 2 above. In the first case there is $p\in\R^m$ such that $p_k\to p$. Then the passage to the limit in
\begin{equation}\label{alt}
\omega_k=\nabla^2\la v_k,h\ra(x_k)\xi_k+\nabla h(x_k)^*p_k,\quad\;(p_k,-\nabla h(x_k)\xi_k)\in\Hat N((h(x_k),v_k);G)
\end{equation}
with taking into account the continuity assumptions and the convergence above, leads to
\begin{equation}\label{alt1}
\omega=\nabla^2\la v,h\ra(\ox)\xi+\nabla h(\ox)^*p,\quad\;(p,-\nabla h(\ox)\xi)\in\Hat N((h(\ox),v);G).
\end{equation}
In the remaining case we have $\|p_k\|\to\infty$ and thus find $p$ such that $p_k/\|p_k\|\to p$ with $\|p\|=1$. Divide now both sides of (\ref{alt}) by $\|p_k\|$ for any large $k$ and take the limit therein as $k\to\infty$. Then we again arrive at (\ref{desc1}). Unifying (\ref{desc1}) and (\ref{alt1}) gives as the description of basic normals to $\O=\gph\partial\ph$ via the following alternative:
\begin{equation}\label{al}
\left[\begin{array}{ll}
(\omega,-\xi)\in N((\ox,\oy);\O)\Longrightarrow\;\exists\,v\in M(\ox,\oy),\;p\in\R^m\;\mbox{ such that}\\\\
\mbox{either: }\;\left\{\begin{array}{ll}
\omega=\nabla^2\la v,h\ra(\ox)\xi+\nabla h(\ox)^*p\\
\mbox{with }\;(p,-\nabla h(\ox)\xi)\in N((h(\ox),v);G)
\end{array}\right.\\\\
\mbox{or: }\;\nabla h(\ox)^*p=0\;\mbox{ with }\;(p,0)\in N(h(\ox),v);G),\;\|p\|=1.
\end{array}\right.
\end{equation}
Remembering the notation introduced in the theorem and the definitions of the constructions used, we see that the ``either" part of (\ref{al}) amounts to the second-order subdifferential inclusion (\ref{am3}) while the ``or" part of (\ref{al}) means the negation of the basic second-order qualification condition (\ref{2qc1}). Thus the assumed fulfillment of (\ref{2qc1}) shows that the ``or" part of (\ref{al}) does not hold, which justifies the validity of the second-order chain rule (\ref{am3}).

Repeating finally the arguments above with taking in to account that the partial counterparts of the first-order chain rule equality (\ref{1-am}) are satisfied due to the results of \cite[Proposition~3.4]{lpr} (see also \cite[Corollary~10.11]{rw} and \cite[Corollary~3.44]{m06a} in more general settings), we get the partial second-order subdifferential chain rules (\ref{am1}) and (\ref{am2}) under the partial second-order qualification condition (\ref{2qc}). $\h$

\begin{Remark}{\bf (second-order chain rules with inclusions).}\label{compar}
{\rm A chain rule in form (\ref{am1}) for the {\em full} second-order subgradient sets of strongly amenable compositions with compatible  parametrization in finite dimensions was derived in \cite{lm} under the second-order subdifferential condition of type (\ref{2qc}) with $\nabla h(\ox,\ow)$ replacing $\nabla_x h(\ox,\ow)$. The proof in \cite{lm} was based on applying a coderivative chain rule to full first-order subdifferential mappings. A similar approach was employed in \cite{m02} and \cite[Theorem~3.74]{m06a} to derive second-order chain rules of type (\ref{am3}) in infinite dimensions under an appropriate infinite-dimensional counterpart of the second-order qualification condition (\ref{2qc1}). Although the results of \cite{m02,m06a} are applied to a more general class of subdifferential regular functions $\th$ in (\ref{am3}), they require a number of additional assumptions in both finite and infinite dimensions. Finally, we mention a second-order chain rule of the inclusion type (\ref{am3}) obtained in \cite[Theorem~3.1]{mo1} for s special kind of strongly amenable compositions with the indicator function $\th=\delta(\cdot;\Theta)$ of a set $\Theta$ in finite dimensions, which does not generally require the fulfillment of the second-order qualification condition (\ref{2qc1}) while imposing instead of a certain {\em calmness} assumption on some auxiliary multifunction. The latter holds, in particular, in the case of polyhedral sets $\Theta$ due to seminal results of \cite{rob1}.}
\end{Remark}

Next we show that the second-order chain rule formula (\ref{am3}), and hence those in (\ref{am1}) and (\ref{am2}), cannot be generally used for {\em precise calculating} the second-order subdifferentials of strongly amenable compositions: the inclusion therein may be {\em strict} even for fairly simple functions $\th$ and $h$ in $\ph=\th\circ h$ without a kind of full rank condition.

\begin{Example}{\bf (strict inclusion in the second-order chain rule formula).}\label{str}
The inclusion in {\rm(\ref{am3})} can be strict even when $h$ is linear while $\th$ is piecewise linear and convex. Moreover, the set on right-hand side of {\rm(\ref{am3})} can be nonempty when the one on the left-hand side is empty.
\end{Example}
{\bf Proof.} Let the functions $h\colon\R^2\to\R^4$  and $\th\colon\R^4\to\R$  be given by
$$
h(x_1,x_2):=(x_1,-x_1,x_2,-x_2)=Ax\;\mbox{ with }\;A:=\left[\begin{array}{c}
\;\;\;1\;\;\;\;\;0\\
-1\;\;\;\;0\\
\;\;\;0\;\;\;\;\;1\\
\;\;\;0\;-1
\end{array}\right],
$$
$$
\th(z_1,z_2,z_3,z_4):=\max\Big\{z_1,z_2,z_3,z_4\Big\}=\sigma_M(z),
$$
where $M:=\{v=(v_1,v_2,v_3,v_4)\in\R^4|\;v_i\ge 0,\;\sum_{i=1}^4 v_i=1\}$ is the unit simplex in $\R^4$, and where $\sigma_\O$ stands for the support function of the set $\O$.  Considering the composition $\ph(x):=\th(h(x))$ on $\R^2$, observe that it can be represented as
\begin{equation}\label{rep}
\ph(x)=\sigma_B(x)\;\mbox{ for }\;B:=\Big\{(y_1,y_2)\in\R^2\Big|\;|y_1|+|y_2|\le 1\Big\}.
\end{equation}
Note that the outer function $\th$ in the strongly amenable (in fact fully amenable) composition $\th\circ h$ is {\em convex piecewise linear} in terminology of \cite{rw}; it can be equivalently described by \cite[Theorem~2.49]{rw} as a function with the polyhedral epigraph.

Using the explicit form (\ref{rep}) of $\ph$ allows us to compute its second-order subdifferential $\partial^2\ph(\ox,\oy)$ with $\ox=(0,0)$ and $\oy=(0,0)$ directly by Definition~\ref{2nd}. Indeed, we get from (\ref{rep}) that $\partial\ph(\ox)=B$, and hence $\oy\in{\rm int}\,\partial\ph(\ox)$. This tells us that
\begin{equation}\label{rep1}
\partial^2\ph(\ox,\oy)(u)=\left\{\begin{array}{ll}
\R^2&\mbox{if }\;u=(0,0),\\
\emp &\mbox{if }\;u\ne(0,0).
\end{array}\right.
\end{equation}
On the other hand, formula (\ref{am3}) reads as the inclusion $\partial^2\ph(\ox,\oy)(u)\subset Q(u)$ with
$$
Q(u):=\disp\bigcup\Big\{A^*\partial^2\th(0,v)(Au)\Big|\;v\in M,\;A^*v=0\Big\}.
$$
Take $\ou=(0,1)$, $\ov=(1/2,1/2,0,0)$ and then check that $A\ou=(0,0,1,-1)$, $\ov\in M$, and $A^*\ov=0$. This ensures the converse inclusion
$$
Q(\ou)\supset A^*\partial^2\th(0,\ov)(A\ou),
$$
which shows that we have $Q(\ou)\ne\emp$ provided that $\partial^2\th(0,\ov)(A\ou)\ne\emp$. To check the latter, recall the representation of the outer function $\th=\sigma_M=\delta_M^*$ via the indicator function $\delta_M=\delta(\cdot;M)$ of $M$. Thus we get the description
$$
\omega\in\partial^2\th(0,\ov)(A\ou)\Longleftrightarrow-A\ou\in\partial^2\delta_M(\ov,0)(-\omega).
$$
Since the set $M$ is a convex polyhedron, an exact formula for $\partial^2\dd_M(0,\ov)$ is available from \cite{dr}.

In order to state this formula, we need to deal with the {\em critical cone} for a convex polyhedron $\O$ at $x\in\O$ with respect to $p\in\partial\dd_\O(x)=N_\O(x)$; this is a polyhedral cone defined by
$$
K(x,p):=\big\{w\in T_\O(\ox)\big|\;w\perp p\big\},
$$
where $T_\O(x)$ is the tangent cone (\ref{2.6}) to $\O$ at $x$. Recall that a {\em closed face} $C$ of a polyhedral cone $K$ is a polyhedral cone of the form
$$
C:=\big\{x\in K|\;x\perp v\}\;\mbox{ for some }\;v\in K^*,
$$
where $K^*$ denotes the polar of the cone $K$. By the proof of \cite[Theorem~2]{dr} (see also \cite[Proposition~4.4]{pr96}) we have the following description of the second-order subdifferential of the indicator function for a convex polyhedron:
\begin{equation}\label{2ind}
w\in\partial^2\dd_\O(x,p)(u)\Longleftrightarrow\left\{\begin{array}{ll}\mbox{there exist closed faces }\;C_1\subset C_2\;\mbox{ of }\; K(x,p)\\
\mbox{with }\;u\in C_1-C_2,\;w\in(C_2-C_1)^*.
\end{array}\right.
\end{equation}
Applying this to our setting with the simplex $\O=M$, we get the critical cone
$$
K=T_M(\ov)\cap 0^\perp=T_M(\ov)=\Big\{(\omega_1,\omega_2,\omega_3,\omega_4)\left|\begin{array}{ll}\omega_1+\omega_2=0\\
\omega_3\ge 0,\;\omega_4\ge 0
\end{array}\right.\Big\}.
$$
It follows from the second-order subdifferential formula (\ref{2ind}) that
$$
-A\ou\in\partial^2\dd_M(\ov,0)(-\omega)\Longleftrightarrow\left\{\begin{array}{ll}\mbox{there exist closed faces }\;C_1\subset C_2\;\mbox{ of }\; K\\
\mbox{with }\;\omega\in C_2-C_1,\;-A\ou\in(C_2-C_1)^*.
\end{array}\right.
$$
Observe that the closed faces of $K$ have the form
$$
\{(\omega_1,\omega_2)|\;\omega_1+\omega_2=0\}\times I\times J,\;\mbox{ where }\;I,J\;\mbox{ can be either }\;\R_+\;\mbox{ or }\;\{0\}.
$$
Denoting by $L$ the subspace in $\{\ldots\}$ of the formula above, we have the following possibilities:
$$
C_1-C_2=L\times(I_1-I_2)\times(J_1-J_2)\;\mbox{ and }\;(C_1-C_2)^*=L^\perp\times(I_1-I_2)^*\times(J_1-J_2)^*,
$$
where $I_1-I_2$ and $J_1-J_2$ can be $\R$, $\R_+$, and $\{0\}$ while, respectively, $(I_1-I_2)^*$ and $(J_1-J_2)^*$ can be $\{0\}$, $\R_-$, and $\R$. Setting now $(I_1-I_2)^*=\R_-$ and $(J_1-J_2)^*=\R$, we get
$$
-A\ou\in\partial^2\dd_M(\ov,0)(-\omega),\;\mbox{ or equilvalently }\;\omega\in\partial^2\th(0,\ov)(A\ou)
$$
whenever $\omega\in L\times\R_+\times\{0\}$. Taking, e.g., $\omega=(0,0,1,0)$ gives us $\partial^2\th(0,\ov)(A\ou)\ne\emp$. Hence the set $Q(\ou)$ on the right-hand side of (\ref{am3}) is nonempty while $\partial^2\ph(\ox,\oy)(\ou)=\emp$ by (\ref{rep1}). $\h$\vspace*{0.05in}

It is not hard to check that the second-order qualification condition (\ref{2qc1}) does not hold in Example~\ref{str}. Thus, besides the emphasis above, this example can be considered as a counterexample to equality in the second-order chain rule (\ref{am3}) with no full rank condition on the derivative and also as an illustration of the possible validity of the inclusion in (\ref{am3}) without the second-order qualification condition (\ref{2qc1}). In the next section we show that the fulfillment of (\ref{2qc1}) yields, for a large class of outer functions $\th$ in amenable compositions $\th\circ h$ including the one in Example~\ref{str}, that the full rank condition must be satisfied, and so the second-order chain rule formula (\ref{am3}) holds in fact as equality.

\section{Analysis of the Basic Second-Order Qualification Condition and Calculating Second-Order Subdifferentials}\sce

This section is mainly devoted to analyzing the basic second-order qualification condition (\ref{2qc1}); our analysis equally works for the partial second-order counterpart (\ref{2qc}).

It is clear that (\ref{2qc1}) holds under the full rank condition
\begin{equation}\label{f-rank}
{\rm rank}\,\nabla h(\ox)=m\Longleftrightarrow\ker\nabla h(\ox)^*=\{0\},
\end{equation}
which ensures the equality in (\ref{am3}) with a unique vector $v$ in (\ref{am3}) and (\ref{2qc1}) by Theorem~\ref{full}. Another setting where condition (\ref{2qc1}) automatically holds is when $\th$ is of class ${\cal C}^{1,1}$ around $\oz=h(\ox)$, which however excludes the case of extended-real-valued outer functions typical in applications of amenable compositions in variational analysis and optimization.

We show in what follows that, for large classes of amenable compositions, the second-order qualification condition (\ref{2qc1}) yields in fact the validity of the full rank condition (\ref{f-rank}) and hence the exact second-order subdifferential chain rule formula
\begin{equation}\label{am4}
\partial^2\ph(\ox,\oy)(u)=\Big(\nabla^2\la\ov,h\ra(\ox)u+\nabla h(\ox)^*\partial^2\th(\oz,\ov)(\nabla h(\ox)u)\Big)\;\mbox{ for all }\;u\in\R^n
\end{equation}
with a unique vector $\ov\in\partial\th(\oz)$ satisfying $\nabla h(\ox)^*v=\oy$. As a by-product of our approach, we efficiently calculate the second-order subdifferential of functions belonging to such classes.

Recall \cite{pr92,rw} that a strongly amenable function $\ph$ is {\em fully amenable} at $\ox$ if the outer function $\th\colon\R^m\to\oR$ in its composite representation $\ph=\th\circ h$ can be chosen as {\em piecewise linear quadratic}. The latter class includes piecewise linear functions discussed in Example~\ref{str} and have their domains as well as their subgradient sets (\ref{2.2}) and (\ref{sin}) to be polyhedral; see \cite[Section10E]{rw}.

We start with a local reduction lemma that describes a general setting where the second-order chain rule (\ref{am3}) holds as equality. Then we show that it is the case for some major classes of fully amenable functions under the second-order qualification condition (\ref{2qc1}).

\begin{Lemma}{\bf (local reduction to full rank condition).}\label{reduc} Let $\ph\colon\R^n\to\oR$ be a strongly amenable composition at $\ox$ represented as $\ph=\th\circ h$ near $\ox$ with $h\colon\R^n\to\R^m$ and $\th\colon\R^m\to\oR$. Denote by $S(z)$ the subspace parallel to the affine hull $\aff\partial\th(z)\subset\R^m$, and let $L$ be any subspace of $\R^m$ satisfying the inclusion
\begin{equation}\label{L}
L\supset S(z)\;\mbox{ for all $z\in\R^m$ sufficiently close to }\;\oz=h(\ox).
\end{equation}
Then we have the exact second-order chain rule formula {\rm(\ref{am4})}.
\end{Lemma}
{\bf Proof.} Let $\dim L=s\le m$. It is easy to see that $s=m$ corresponds to the full rank condition on $\nabla h(\oz)$, and thus the second-order chain rule (\ref{am4}) follows from Theorem~\ref{full}. Suppose now that $s<m$, and let $A$ be the matrix of a linear isometry from $\R^m$ into $\R^s\times\R^{m-s}$ under which $AL=\R^s\times\{0\}$. Denoting $P:=Ah$ and $\vartheta:=\th A^{-1}$ gives us the representation $\ph=\vartheta\circ P$. Hence the initial framework of the lemma can be reduced to one in which we have, in terms of $P(x)=(p_1(x),\ldots,p_m(x))$ and $z=Ax$, the implication
$$
\left.\begin{array}{ll}
z\;\mbox{ sufficiently close to }\;\oz\\
v=(v_1,\ldots,v_m)\in\partial\vartheta(z)
\end{array}\right\}\Longrightarrow v_{s+1}=0,\ldots,v_m=0.
$$
This means that in analyzing $\partial\ph$ locally via $\vartheta$ and $P$ it is possible to pass with no loss of generality to the ``submapping"
$$
P_0\colon x\mapsto\big(p_1(x),\ldots,p_s(x)\big),
$$
since only $p_1,\ldots,p_s$ are active locally while $p_{s+1},\ldots,p_m$ do not matter in the implication
$$
y\in\partial\ph(x)\Longrightarrow\exists\,v\in\partial\vartheta(P(x))\;\mbox{ such that }\;\nabla P(x)^*v=y.
$$
It suffices therefore to impose a full rank $(=s$) condition on $P_0$ and invoke the second-order chain rule in equality form (\ref{am4}) from Theorem~\ref{full}. The result of the lemma can be then translated back to the original context of $\ph=\th\circ h$ through $A^{-1}$. $\h$\vspace*{0.05in}

Observe that the full rank condition on $P_0$ comes out in the context of Lemma~\ref{reduc} from
$$
L\cap\ker\nabla h(x)^*=\{0\}.
$$
We show in what follows that a subspace $L$ in (\ref{L}) is provided by the second-order qualification condition (\ref{2qc1}) when $\ph=\th\circ h$ belongs to major classes of fully amenable compositions.\vspace*{0.05in}

To proceed in this direction, let us first calculate the second-order subdifferential $\partial^2\th(\oz,\ov)(0)$ from the left-hand side of (\ref{2qc1}) for general convex piecewise linear-quadratic functions $\th$, which is certainly of its own interest.

\begin{Theorem}{\bf (calculating the second-order subdifferential of convex linear-quadratic functions).}\label{lq}
Let $\ph=\th\circ h$ be a fully amenable composition at $\ox$, let $M\colon\R^n\times\R^n\tto\R^m$ be a set-valued mapping defined in {\rm(\ref{M})}, and let $S(z)$ be a subspace of $\R^m$ parallel to the affine hull $\aff\partial\th(z)$ for $z$ near $\oz=h(\ox)$. Then for any sufficiently small neighborhood $O$ of $\oz$ there is a finite union of the subspaces $S(z)$ such that
\begin{equation}\label{un}
\partial^2\th(\oz,\ov)(0)=\bigcup_{z\in O}S(z)\;\mbox{ whenever }\;\ov\in M(\ox,\oy).
\end{equation}
\end{Theorem}
{\bf Proof.} As mentioned in the proof of Theorem~\ref{am}, the mapping $M$ from (\ref{M}) is closed-graph around $(\ox,\oy,\oz)$ and uniformly bounded around $(\ox,\oy)$ for strongly amenable compositions. Since $\th\circ h$ is fully amenable, the set $M(x,y)$ is also polyhedral for all $(x,y)$ sufficiently close to $(\ox,\oy)$.

Fix any $\ov\in M(\ox,\oy)$. Since $\th$ is piecewise linear-quadratic, its graph $G:=\gph\th$ is {\em piecewise polyhedral}, i.e., it is the union of finitely many polyhedral sets in $\R^m$. Using this and taking formulas (\ref{2.4}) and (\ref{2.5}) into account, we find a neighborhood  $W$ of $(\oz,\ov)$ such that
\begin{equation}\label{poly}
N_G(\oz,\ov)=\bigcup\Big\{\Hat N_G(z,v)\Big|\;(z,v)\in G\cap W\Big\}=\bigcup\Big\{T_G(z,v)^*\Big|\;(z,v)\in G\cap W\Big\},
\end{equation}
where only {\em finitely many} cones (all of them are polyhedral) occur in the unions. Therefore
\begin{equation}\label{poly1}
w\in\partial^2\th(\oz,\ov)(0)\Longleftrightarrow\exists\,(z,v)\in G\cap W\;\mbox{ with }\;(w,0)\in T_G(z,v)^*.
\end{equation}
On the other hand, we have $T_G(z,v)=\gph(D\partial\th)(z,v)$ by definition of the {\em graphical derivative} $D$ of a set-valued mapping, and furthermore
$$
(D\partial\th)(z,v)=\partial\Big(\hbox{$1\over
2$}d^2\th(z,v)\Big)
$$
via the {\em second subderivative} of the function $\th$ under consideration; see \cite[Theorem~13.40 and Proposition~13.32]{rw} for moire details. Hence it ensures that
$$
\dom(D\partial\th)(z,v)=\dom d^2\th(z,v)=N_{\partial\th(z)}(v)
$$
by \cite[Theorem~13.14]{rw}. Employing now (\ref{poly}) and (\ref{poly1}) gives us the representations
\begin{equation}\label{poly2}
\begin{array}{ll}
\partial^2\th(\oz,\ov)(0)&=\disp\bigcup_{(z,v)\in G\cap W}\Big[\dom(D\partial\th)(z,v)\Big]^*=\bigcup_{(z,v)\in G\cap W}\Big[N_{\partial\th(u)}(v)\Big]^*\\
&=\disp\bigcup_{(z,v)\in G\cap W}T_{\partial\th(z)}(v),
\end{array}
\end{equation}
where only finitely many sets are taken in the unions. Pick $v\in\partial\th(z)$ and find, by the polyhedrality of the subgradient sets $\partial\th(z)$ and the construction of the subspaces $S(z)$, a vector $v'\in\ri\partial\phi(z)$ arbitrary close to $v'$ and get for all such $v'$ the relationships
$$
T_{\partial\th(z)}(v')=S(z)\supset T_{\partial\th(z)}(v),
$$
which imply by (\ref{poly2}) the equality
$$
\partial^2\th(\oz,\ov)(0)=\bigcup S(z),
$$
where the finite union of {\em subspaces} are taken over $z$ such that $(z,v)\in G\cap W$ for some $v$.

So far we focused our analysis on a particular point $\ov\in M(\ox,\oy)$ and an associated neighborhood $W$ of $(\oz,\ov)$. Since the mapping $M$ is of closed graph and locally bounded, the set $M(\ox,\oy)$ can be covered by {\em finitely many} of such neighborhoods. This allows us to obtain (\ref{un}) and complete the proof of the theorem. $\h$\vspace*{0.05in}

Note that representation (\ref{un}) held for outer functions of any fully amenable composition is not generally sufficient for applying Lemma~\ref{reduc} and deducing thus the second-order chain rule as equality (\ref{am4}) from the second-order qualification condition (\ref{2qc1}). To proceed in this direction, we need to get just {\em one subspace} in representation (\ref{un}), which serves all $\ov\in M(\ox,\oy)$ therein. The next result shows that it can be done in the case of {\em piecewise linear} outer functions in fully amenable compositions $\ph=\th\circ h$.

\begin{Theorem}{\bf (second-order chain rule for fully amenable compositions with piecewise linear outer functions).}\label{pl} Let $\ph=\th\circ h$ be a fully amenable composition at $\ox$, where $\th\colon\R^m\to\oR$ is $($convex$)$ piecewise linear. Assume that the second-order qualification condition {\rm(\ref{2qc1})} is satisfied. Then we have the exact second-order chain rule formula {\rm(\ref{am4})}.
\end{Theorem}
{\bf Proof.} It is not hard to check that for convex piecewise linear functions $\ph$ we have the inclusion $\partial\th(z)\subset\partial\th(\oz)$ for any neighborhood $O$ of $\oz$ sufficiently small. This implies that $S(z)\subset S(\oz)$ whenever $z\in O$ and hence ensures the equality
\begin{equation}\label{2}
\partial^2\th(\oz,\ov)(0)=S(\oz)
\end{equation}
by formula (\ref{un}) from Theorem~\ref{lq}. Employing the latter subdifferential representation in the second-order qualification condition (\ref{2qc1}) gives us
\begin{equation}\label{2nd2}
S(z)\cap\ker\nabla h(\ox)^*=\{0\}\;\mbox{ for all }\;z\in O.
\end{equation}
Recall that the subspace $S(z)$ consists of all vectors $\lm(v'-v)$ such that $\lm\in\R$ and $\lm,\lm'\in\partial\th(z)$. Hence the second-order qualification condition (\ref{2nd2}) is equivalent to the following: there exist neighborhoods $U$ of $\oz$ and $V$ of $\oy$ such that
$$
\Big[x\in U,\;y\in V,\;v,v'\in M(x,y)\Big]\Longrightarrow v=v'.
$$
On other words, the latter means that the mapping $M$ from (\ref{M}) is {\em single-valued} on the subdifferential graph $\gph\partial\ph$ around $(\ox,\oy)$ in the case under consideration; thus it is continuous as well. The result of the theorem follows now from Lemma~\ref{reduc}. $\h$\vspace*{0.05in}

Next we consider a major subclass of piecewise linear-quadratic outer functions in fully amenable compositions given by
\begin{equation}\label{maj}
\th(z):=\sup_{v\in C}\Big\{\la v,z\ra-\hbox{$1\over 2$}\la v,Qv\ra\Big\},
\end{equation}
where $C\subset\R^m$ is a nonempty polyhedral set, and where $Q\in\R^{m\times m}$ is a symmetric positive-semidefinite matrix. Functions of this class are useful in many aspects if variational analysis and optimization; in particular, as penalty expressions in composite formats of optimization; see, e.g., \cite{rw} and the references therein. By definition (\ref{maj}) we see that $\th$ in (\ref{maj}) is proper, convex, and piecewise linear-quadratic (piecewise linear when $Q=0$) with the conjugate representation
\begin{equation}\label{conj}
\th(z)=(\dd_C+j_Q)^*(z)\;\mbox{ for }\;j_Q(v):=\hbox{$1\over 2$}\la v,Qv\ra.
\end{equation}
In the following theorem we calculate the second-order subdifferential of functions $\th$ from (\ref{maj}) and justifies the equality second-order chain rule formula (\ref{am4}) for fully amenable compositions $\ph=\th\circ h$ with outer functions of this type.

\begin{Theorem}{\bf (second-order calculus rule for a major subclass of fully amenable compositions).}\label{spec}
Let $\ph=\th\circ h$ be a fully amenable composition at $\ox$ with $\th\colon\R^m\to\oR$ of class {\rm(\ref{maj})} satisfying the assumptions made above. Suppose also that the second-order qualification condition {\rm(\ref{2qc1})} is satisfied. Then we have the exact second-order chain rule formula {\rm(\ref{am4})}.
\end{Theorem}
{\bf Proof.} It follows from (\ref{conj}) that the conjugate to $\th$ is $\th^*=\dd_C+j_Q$. Observe also that
\begin{equation}\label{con}
w\in\partial^2\th(\oz,\ov)(u)\Longleftrightarrow-u\in\partial^2\th^*(\ov,\oz)(-w).
\end{equation}
Furthermore, we have by the calculations in \cite[Example~11.18]{rw} that
\begin{equation}\label{con1}
\partial\th^*(v)=N_C(v)+Qv,\quad v\in C,
\end{equation}
and hence $\oz\in\partial\th^*(\ov)\Longleftrightarrow\oz-Q\ov\in N_C(\ov)$. Using this and definition (\ref{2nd1}) of the second-order subdifferential and then applying the coderivative sum rule \cite[Theorem~1.62]{m06a} to (\ref{con1}) give us
$$
\partial^2\th^*(\ov,\ou)(-w)+D^* N_C(\ov,\oz-Q\ov)(-w)+\nabla j_Q(-w).
$$
Since $\nabla j_Q(-w)=Qw$, the latter implies in turn that
\begin{equation}\label{con2}
-u\in\partial^2\th^*(\ov,\oz)(-w)\Longleftrightarrow Qw-u\in\partial^2\dd_C(\ov,\oz-Q\ov)(-w).
\end{equation}
Employing (\ref{con}, (\ref{con2}) and proceeding similarly to the consideration in Example~\ref{str} above, we get from (\ref{con2}) the {\em exact formula for calculating the second-order subdifferential}:
\begin{equation}\label{con3}
u\in\partial^2\th(\oz,\ov)(w)\Longleftrightarrow\left\{\begin{array}{ll}
\exists\,\mbox{closed faces }\;K_1\supset K_2\;\mbox{ of }\;K\\
\mbox{with }\;w\in K_1-K_2,\;Qw-u\in(K_1-K_2)^*,
\end{array}\right.
\end{equation}
where $K$ is the critical cone for $C$ at $\ov$ with respect to $\oz-Q\ov$ given by
$$
K=T_C(\ov)\cap(\oz-Q\ov)^\perp.
$$
The positive-semidefiniteness of the matrix $Q$ yields that
$$
0\ge\la w,Qw\ra\Longleftrightarrow Qw=0\Longleftrightarrow z\in\ker Q,
$$
which allows us to deduce from (\ref{con3}) that
\begin{eqnarray*}\begin{array}{ll}
w\in\partial^2\th(\oz,\ov)(0)&\Longleftrightarrow\;\exists\,K_1\supset K_2\;\mbox{ with }\;w\in K_1-K_2,\;Qw\in(K_1-K_2)^*\\
&\Longleftrightarrow\;\exists\,K_1\supset K_2\;\mbox{ with }\;w\in(\ker Q)\cap(K_1-K_2)\\
&\Longleftrightarrow\;w\in(\ker Q)\cap(K-K).
\end{array}
\end{eqnarray*}
Thus the set $\partial^2\th(\oz,\ov)(0)$ is a {\em subspace} in $\R^m$. Substituting it finally into the second-order subdifferential condition (\ref{2qc1}) and taking into account that  $\partial\th(z)=(N_C+Q)^{-1}(z)$ by \cite[Example~11.18]{rw}, we arrive at the second-order equality chain rule (\ref{am4}) similarly to the proof of Theorem~\ref{pl} based on the application of Lemma~\ref{reduc}. $\h$

\section{Applications to Tilt Stability in Nonlinear and Extended Nonlinear Programming}\sce

The second-order chain rules and subdifferential calculations obtained in Sections~3 and 4 are undoubtedly useful in any settings where the second-order subdifferential (\ref{2nd1}) and its partial counterparts are involved; see the discussions and references in Section~1. In this section we confine ourselves to the usage of second-order chain rules for deriving {\em full characterizations} of {\em tilt-stable local minimizers} in some important classes of constrained optimization problems. It requires applying {\em equality-type} formulas of the second-order subdifferential calculus.

The notion of title-stable minimizers was introduced by Poliquin and Rockafellar \cite{pr98} in order to characterize strong manifestations of optimality that support computational work via the study of how local optimal solutions react to shifts (tilt perturbations) of the data. Recall that a point $\ox$ is a {\em tilt-stable local minimizer}{ of the function $\ph\colon\R^n\to\oR$ finite at $\ox$ if there is $\gg>0$ such that mapping
$$
M\colon y\mapsto\disp{\rm argmin}\Big\{\ph(x)-\ph(\ox)-\la y,x-\ox\ra\Big|\;\|x-\ox\|\le\gg\Big\}
$$
is single-valued and Lipschitz continuous on some neighborhood of $y=0$ with $M(0)=\ox$.

It is proved in \cite[Theorem~1.3]{pr98} that for $\ph\colon\R^n\to\oR$ having $0\in\partial\ph(\ox)$ and such that $\ph$ is both prox-regular and subdifferentially continuous at $\ox$ for $\oy=0$, the point $\ox$ is a tilt-stable local minimizer of $\ph$ if and only if the second-order subdifferential mapping $\partial^2\ph(\ox,0)\colon\R^n\tto\R^n$ is
{\em positive-definite} in the sense that
\begin{equation}\label{pd}
\la w,u\ra>0\;\mbox{ whenever }\;w\in\partial^2\ph(\ox,0)(u)\;\mbox{ with }\;u\ne 0.
\end{equation}
The aforementioned properties of prox-regularity and subdifferential continuity introduced in \cite{pr96} (see also \cite[Definitions~13.27 and 13.28]{rw}) hold for broad classes of ``nice" functions encountered in variational analysis and optimization. In particular, both properties are satisfied, at all points of a neighborhood of $\ox$ for any function strongly amenable at $\ox$; see \cite[Proposition~13.32]{rw}.

Our subsequent goal is to extend the characterization of tilt-stable local minimizers from \cite{pr98} to favorable classes of {\em constrained optimization} problems.
To proceed, we use the following {\em composite format} of optimization known as {\em extended nonlinear programming} (ENLP); see \cite{r00,rw}:
\begin{equation}\label{copt}
\mbox{minimize }\;\ph(x):=\ph_0(x)+\th\big(\ph_1(x),\ldots,\ph_m(x)\big)=\ph_0(x)+(\th\circ\Phi)(x)\;\mbox{ over }\;x\in\R^n,
\end{equation}
where $\th\colon\R^m\to\oR$ is an extended-real-valued function, and where $\Phi(x):=(\ph_1(x),\ldots\ph_m(x))$ is a mapping from $\R^n$ to $\R^m$. Written in the unconstrained format, problem (\ref{copt}) is actually a problem of constrained optimization with the set of feasible solutions given by
$$
X:=\{x\in\R^n|\;(\ph_1(x),\ldots,\ph_m(x))\in Z\}\;\mbox{ for }\;Z:=\{z\in\R^m|\;\th(z)<\infty\}.
$$
In other words, problem (\ref{copt}) can be equivalently represented in the form
\begin{equation}\label{copt1}
\mbox{minimize }\;\ph_0(x)+\dd_Z(\Phi(x))\;\mbox{ over }\;x\in\R^m\;\mbox{ with }\;Z=\dom\th
\end{equation}
via the indicator function of the feasible set. As argued in \cite{r00}, the composite format (\ref{copt}), or (\ref{copt1}), is a convenient  general framework from both theoretical and computational viewpoints to accommodate a variety of particular models in constrained optimization. Note that the conventional problem of nonlinear programming with $s$ inequality constraints and $m-s$ equality constraints can be written in form (\ref{copt1}), where $Z=\R^s_-\times\{0\}^{m-s}$.

Our first result provides a complete second-order characterization of tilt-stable minimizers $\ox$ for a general class of problems (\ref{copt}) under full rank of the Jacobian matrix $\nabla\Phi(\ox)$.

\begin{Theorem} {\bf (characterization of tilt-stable minimizers for constrained problems with full rank condition).}\label{fr} Let $\ox\in X$ be a feasible solution to {\rm(\ref{copt})} such that $\ph_0$ and $\Phi$ are smooth around $\ox$ with their derivatives strictly differentiable at $\ox$, that ${\rm rank}\,\nabla\Phi(\ox)=m$, and that $\th$ is prox-regular and subdifferentially continuous at $\oz:=\Phi(\ox)$ for the $($unique$)$ vector $\ov\in\R^m$ satisfying the relationships
\begin{equation}\label{fr1}
\ov\in\partial\th(\oz)\;\mbox{ and }\;\nabla\Phi(\ox)^*\ov=-\nabla\ph_0(\ox).
\end{equation}
Then $\ox$ with $-\nabla\ph_0(\ox)\in\nabla\Phi(\ox)^*\partial\th(\oz)$ is a tilt-stable local minimizer of {\rm(\ref{copt})} if and only if the mapping $T\colon\R^n\tto\R^n$ given by
\begin{equation}\label{T}
T(u):=\nabla^2\la\ov,\Phi\ra(\ox)(u)+\nabla\Phi(\ox)^*\partial^2\th(\oz,\ov)(\nabla\Phi(\ox)u),\quad u\in\R^n,
\end{equation}
is positive-definite in the sense of {\rm(\ref{pd})}.
\end{Theorem}
{\bf Proof.} Since $\ph_0$ and $\Phi$ are smooth around $\ox$ and $\nabla\Phi(\ox)$ has full rank $m$, it follows from the first-order subdifferential sum and chain rules of \cite[Corollary~10.9 and Exercise~10.7]{rw} that
$$
0\in\partial\ph(\ox)\Longleftrightarrow-\nabla\ph_0(\ox)\in\nabla\Phi(\ox)^*\partial\th(\oz)
$$
for $\ph$ in (\ref{copt}). Furthermore, these rules and the definitions of prox-regularity and subdifferential continuity in \cite{rw} imply that the latter properties of $\ph$ at $\ox$ for $0\in\partial\ph(\ox)$ are equivalent to the corresponding properties of $\th$ at $\oz$ for $\ov$ satisfying (\ref{fr1}).

It remains to check therefore that the positive-definiteness (\ref{pd}) of $\partial^2\ph(\ox,0)$ is equivalent to that of $T$ in (\ref{T}). We show in fact that $\partial^2\ph(\ox,0)(u)=T(u)$ for all $u\in\R^n$. Indeed, using the second-order chain rule from \cite[Proposition~1.121]{m06a} in (\ref{copt}) gives us
\begin{equation}\label{fr0}
\partial^2\ph(\ox,0)(u)=\nabla^2\ph_0(\ox)u+\partial^2(\th\circ\Phi)(\ox,-\nabla\ph_0(\ox))(u),\quad u\in\R^n.
\end{equation}
To complete the proof of the theorem, we finally apply the exact second-order chain rule from Theorem~\ref{full} to the composition $\th\circ\Phi$ in the latter equality.
$\h$\vspace*{0.05in}

Next we address the conventional model of {\em nonlinear programming} (NLP) with smooth data:
\begin{equation}\label{nlp}
\mbox{minimize }\;\ph_0(x)\;\mbox{ subject to }\;\ph_i(x)=\left\{\begin{array}{ll}
\le 0 &\mbox{for }\;i=1,\ldots,s,\\
=0&\mbox{for }\;i=s+1,\ldots,m.
\end{array}\right.
\end{equation}
As mentioned above, problem (\ref{nlp}) can be written in form (\ref{copt1}) with $Z=\R^s_-\times\{0\}^{m-s}$. For this problem, the full rank condition of Theorem~\ref{fr} corresponds to: the gradients
\begin{equation}\label{li}
\nabla\ph_1(\ox),\ldots,\nabla\ph_m(\ox)\;\mbox{ are linearly independent}.
\end{equation}
Actually, since our analysis is local, we can drop in what follows any inactive inequality constraints from the picture and thus reduce with no loss of generality to having
\begin{equation}\label{ac}
\ph_i(\ox)=0\;\mbox{ for all }\;i=1,\ldots,m.
\end{equation}
The full rank condition (\ref{li}) in case (\ref{ac}) is then the classical (LICQ): {\em the active constraint gradients at $\ox$ are linearly independent}.

To proceed further, consider the Lagrangian function in (\ref{nlp}) given by
$$
L(x,\lm):=\ph_0(x)+\sum_{i=1}^m\lm_i\ph_i(x)\;\mbox{ with }\;\lm=(\lm_1,\ldots,\lm_m)\in\R^m
$$
and remember that, for any local optimal solution $\ox$ to (\ref{nlp}), the LICQ at $\ox$ ensures the existence of a {\em unique} multiplier vector $\bar\lm=(\bar\lm_1,\ldots,\bar\lm_m)\in\R^s_+\times\R^{m-s}$ such that
\begin{equation}\label{lm}
\nabla_x L(\ox,\bar\lm)=\nabla\ph_0(\ox)+\sum_{i=1}^m\bar\lm_i\nabla\ph_i(\ox)=0.
\end{equation}
Recall that the {\em strong second-order optimality condition} (SSOC) holds at $\ox$ if
\begin{equation}\label{ssoc}
\la u,\nabla^2_{xx}L(\ox,\bar\lm)u\ra>0\;\mbox{ for all }\;0\ne u\in S,
\end{equation}
where the subspace $S\subset\R^n$ is given by
$$
S:=\{u\in\R^n|\;\la\nabla\ph_i(\ox),u\ra=0\;\mbox{ whenever }\;i=1,\ldots,m\}.
$$
Note that (\ref{ssoc}) is also known as the ``strong second-order sufficient condition" for local optimality. The following theorem shows that, in the setting under consideration, the SSOC is {\em necessary and sufficient} for the tilt stability of local minimizers.

\begin{Theorem}{\bf (characterization of tilt-stable local minimizers for NLP).}\label{t-nlp} Let $\ox$ be a feasible solution to {\rm(\ref{nlp})} such that all the functions $\ph_i$ for $i=0,\ldots,m$ are smooth around $\ox$ with their derivatives strictly differentiable at $\ox$ and that the LICQ is satisfied at this point. Then we have the following assertions:

{\bf (i)} If $\ox$ is a tilt-stable local minimizer of {\rm(\ref{nlp})}, then SSOC {\rm(\ref{ssoc})} holds at $\ox$ with the unique multiplier vector $\bar\lm\in\R^s_+\times\R^{m-s}$ satisfying {\rm(\ref{lm})}.

{\bf (ii)} Conversely, the validity of SSOC at $\ox$ with $\bar\lm\in\R^s_+\times\R^{m-s}$ satisfying {\rm(\ref{lm})} implies that $\ox$ is a tilt-stable local minimizer of {\rm(\ref{nlp})}.
\end{Theorem}
{\bf Proof.} As mentioned above, the LICQ corresponds to the full rank condition of Theorem~\ref{fr} in the setting under consideration. The prox-regularity and subdifferential continuity of the indicator function $\th=\dd_Z$ with $Z=\R^s_-\times\{0\}^{m-s}$ follow from its convexity \cite[Example~13.30]{rw}. Let us next represent the mapping $T$ in (\ref{T}) via the initial data of problem (\ref{nlp}). It is easy to see that $T(u)$ reduces in this case to
$$
T(u)=\nabla^2_{xx}L(\ox,\bar\lm)u+\nabla\Phi(\ox)^*\partial^2\dd_Z(0,\bar\lm)(\nabla\Phi(\ox)u),
$$
with $\Phi=(\ph_1,\ldots,\ph_m)$ and $Z=\R^s_-\times\{0\}^{m-s}$, provided that the first-order condition (\ref{lm}) is satisfied, which is of course the case when $\ox$ is a local minimizer of (\ref{nlp}). Thus the positive-definiteness of $T(u)$ amounts to
\begin{equation}\label{lm1}
u\ne 0,\;w\in\partial^2\dd_Z(0,\bar\lm)(u)\Longrightarrow\la u,\nabla^2_{xx}L(\ox,\bar\lm)u\ra+\la w,\nabla\Phi(\ox)u\ra>0.
\end{equation}
To proceed, we calculate the second-order subdifferential $\partial^2\dd_Z(0,\bar\lm)$ in (\ref{lm1}) by using formula (\ref{2ind}) presented and discussed in Example~\ref{str}. Observe that the critical cone in this situation is $K=Z\cap\bar\lm^{\perp}$. It follows directly from (\ref{2ind}) that
$$
w\in\partial^2\dd_Z(0,\bar\lm)(u)\Longleftrightarrow\left\{\begin{array}{ll}\mbox{there exist closed faces }\;K_1\subset K_2\;\mbox{ of }\;K\\
\mbox{with }-\nabla\Phi(\ox)u\in K_1-K_2,\;w\in(K_2-K_1)^*.
\end{array}\right.
$$
The latter implies in turn that
$$
\min_{w\in\partial^2\dd_Z(0,\bar\lm)(u)}\la w,\nabla\Phi(\ox)u\ra=0\;\mbox{ for all }\;u\in\dom\partial^2\dd_Z(0,\bar\lm)
$$
with the subdifferential domain representation
$$
\dom\partial^2\dd_Z(0,\bar\lm)=\bigcup\Big\{(K_1-K_2)\Big|\;K_1\subset K_2\;\mbox{ closed faces of }\;K\Big\}=K-K.
$$
Substituting this into (\ref{lm1}) and taking into account the forms of the critical cone $K$ as well as the subspace $S$ in (\ref{ssoc}), we conclude that the positive-definiteness condition (\ref{lm1}) is equivalent to the strong second-order optimality condition (\ref{ssoc}) provided that (\ref{lm}) holds at $\ox$.

Having this in hand and using Theorem~\ref{fr} allow us to justify both assertions in (i) and (ii). Indeed, since every tilt-stable local minimizers $\ox$ is a standard local minimizer, it satisfies the first-order necessary optimality condition (\ref{lm}) under the assumed LICQ at $\ox$. Thus  the SSOC holds at this point, which is the assertion in (i). The validity of the converse assertion (ii) follows from the equivalence between (\ref{lm1}) and (\ref{ssoc}) under (\ref{lm}) proved above. $\h$\vspace*{0.05in}

The obtained characterization of tilt-stable minimizers for NLP leads us to comparing this notion with the classical Robinson's notion of {\em strong regularity} \cite{rob} of parameterized variational inequalities associated with the KKT conditions for NLP (\ref{nlp}). Complete characterizations of strong regularity for NLP are derived in \cite{dr}; see also the references therein.

\begin{Corollary}{\bf (comparing tilt-stability and strong regularity).}\label{sreg}
Under the assumptions of Theorem~{\rm\ref{t-nlp}}, the tilt-stability of local minimizers for {\rm(\ref{nlp})} is equivalent to the strong regularity of the variational inequality associated with the KKT conditions for {\rm(\ref{nlp})}.
\end{Corollary}
{\bf Proof}. It follows directly from Theorem~\ref{t-nlp} and the characterization of strong regularity obtained in \cite[Theorem~5 and Theorem~6]{dr}. $\h$\vspace*{0.05in}

It is not hard to check that the strong regularity of the KKT system directly implies the LICQ at the corresponding solution of (\ref{nlp}). On the other hand, the LICQ requirement arising from the full rank condition of Theorem~\ref{fr} is essential for the SSOC characterization of tilt-stability in Theorem~\ref{t-nlp}. Furthermore, even imposing the seemingly less restrictive second-order qualification condition (\ref{2qc1}) needed for deriving the second-order chain rule {\em unavoidably} leads us the the LICQ requirement for NLP, since the latter class is represented via fully amenable compositions with piecewise linear outer functions $\th$ in the composite format (\ref{copt}). This follows from the results of Section~4 and is reflected in the next theorem.

\begin{Theorem}{\bf (characterizing tilt-stable minimizers for constrained problems described by fully amenable compositions).}\label{t-stab}  Let $\ox$ be a feasible solution to {\rm(\ref{copt})} such that $\ph_0$ smooth around $\ox$ with the strictly differentiable derivative at $\ox$ and that the composition $\th\circ\Phi$ is fully amenable at $\ox$ with the outer function $\th\colon\R^m\to\oR$ of the following types:
\begin{itemize}

\item either $\th$ is piecewise linear,

\item or $\th$ is of class {\rm(\ref{maj})} with a nonempty polyhedral set $C\subset\R^m$ and a symmetric positive-semidefinite matrix $Q\in\R^m\times\R^m$.
\end{itemize}
Assume further that the second-order qualification condition {\rm(\ref{2qc1})} holds at $\ox$, where $v=-\ov$ is the unique vector satisfying {\rm(\ref{fr1})} with $\oz=\Phi(\ox)$.
Then $\ox$ with $-\nabla\ph_0(\ox)\in\nabla\Phi(\ox)^*\partial\th(\oz)$ is a tilt-stable local minimizer of {\rm(\ref{copt})} if and only if the mapping $T\colon\R^n\tto\R^n$ defined in {\rm(\ref{T})} is positive-definite in the sense of {\rm(\ref{pd})}, where the second-order subdifferential $\partial^2\th(\oz,\ov)$ is calculated by formulas  {\rm(\ref{2})} and {\rm(\ref{con3})}, respectively.
\end{Theorem}
{\bf Proof.}  Observe first that in both cases under consideration the composition $\th\circ\Phi$ is prox-regular and subdifferentially continuous at any point $x$ around $\ox$ by \cite[Proposition~13.32]{rw}; hence the same holds for the function $\ph$ from (\ref{copt}). It follows from Theorems~\ref{pl} and \ref{spec} that, under the validity of the second-order qualification condition (\ref{2qc1}), we have the unique vector $\ov$ satisfying (\ref{fr1}) and the second-order chain rule
\begin{equation}\label{t-stab1}
\partial^2(\th\circ\Phi)(\ox,-\nabla\ph_0(\ox))(u)=\Big(\nabla^2\la\ov,\Phi\ra(\ox)u+\nabla\Phi(\ox)^*\partial^2\th(\oz,\ov)(\nabla\Phi(\ox)u)\Big)
\end{equation}
for all $u\in\R^n$ when $\th$ belongs to one of the classes considered in this theorem. Substituting further (\ref{t-stab1}) into formula (\ref{fr0}) due to the
the second-order sum rule from \cite[Proposition~1.121]{m06a} allows us to justify that
$$
\partial^2\ph(\ox,0)=T(u)\;\mbox{ whenever }\;u\in\R^n,
$$
and thus the positive-definiteness of the mapping $T$ from (\ref{T}) fully characterizes the tilt-stability of the local minimizer $\ox$ of (\ref{comp}) in both cases of $\th$ under consideration with the formulas for calculating of $\partial^2\th(\oz,\ov)$ derived in the proofs of Theorems~\ref{pl} and \ref{spec}, respectively. $\h$\vspace*{0.05in}

We conclude the paper with the following two final remarks.

\begin{Remark}{\bf (sufficient conditions for tilt-stable local minimizers).}\label{suf}
{\rm The second-order chain rule (\ref{am3}) of the {\em inclusion} type derived in Theorem~\ref{am} for {\em strongly amenable} compositions and the second-order sum rule inclusions obtained in \cite[Theorem~3.73]{m06a} allow us to establish general {\em sufficient} conditions for tilt-stable local minimizers in large classes of constrained optimization problems written in the composite format (\ref{copt}). Indeed if, in addition to the hypotheses of Theorem~\ref{am} for the composition $\th\circ\Phi$ in (\ref{copt}), we assume that the function $\ph_0$ is, e.g., of class ${\cal C}^{1,1}$ around $\ox$, then we have by the second-order sum rule from \cite[Theorem~3.73(i)]{m06a} and the chain rule of Theorem~\ref{am} the fulfillment of the inclusion
\begin{equation}\label{inc}
\partial^2\ph(\ox,0)(u)\subset T(u),\quad u\in\R^n,
\end{equation}
for $\ph$ from (\ref{copt}) and $T$ from (\ref{T}). The prox-regularity and subdifferential continuity of such functions $\ph$ follow, under the assumptions made, from \cite[Proposition~13.32 and Proposition~13.34]{rw} and first-order subdifferential calculus rules. Thus inclusion (\ref{inc}) ensures that the positive-definiteness of $T$ implies the one of $\partial^2\ph(\ox,0)$, and the former is therefore a sufficient condition for the tilt-stability of local minimizers of (\ref{copt}).}
\end{Remark}

\begin{Remark}{\bf (full stability of local minimizers).}\label{fs}
{\rm Developing the concept of tilt stability, Levy, Poliquin and Rockafellar \cite{lpr} introduced the notion of {\em fully stable} local minimizers of general optimization problems of the type
\begin{equation}\label{fs1}
\mbox{minimize }\;\ph(x,u)-\la v,x\ra\;\mbox{ over }\;x\in\R^n
\end{equation}
with respect to both ``basic" perturbations $u$ and ``tilt" perturbations $v$. The main result of that paper \cite[Theorem~2.3]{lpr} establishes a {\em complete characterization} of fully stable local minimizers of (\ref{fs1}) via the positive-definiteness of the {\em extended partial second-order subdifferential} (\ref{par2}) of $\ph$. Similarly to the results of this section for tilt-stable minimizers of constraint optimization problems written in the composite format (\ref{copt}), we can derive characterizations as well as sufficient conditions for fully stable minimizers of (\ref{fs1}) based on \cite[Theorem~2.3]{lpr} and the second-order chain rules for the partial second-order counterpart (\ref{par2}) obtained above. Our ongoing research project is to comprehensively elaborate these developments on full stability in constrained optimization and its applications.}
\end{Remark}

\begin{Remark}{\bf (tilt stability and partial smoothness).} \label{lz} {\rm After completing this paper, we became aware of the concurrent work by Lewis and Zhang \cite{lz} related to second-order subdifferentials (generalized Hessians) and tilt stability. The main results of \cite{lz} provide calculations of the basic second-order construction from Definition~\ref{2nd} for ${\cal C}^2$-partly smooth functions on ${\cal C}^2$-smooth manifolds and then characterize tilt stability in such settings via strong criticality and local quadratic growth.}
\end{Remark}


\small
\begin{thebibliography}{10}

\bibitem{bms} L. BAN, B. S. MORDUKHOVICH and W. SONG, {\em Lipschitzian stability of parametric variational inequalities over generalized polyhedra in Banach spaces}, Nonlinear Anal.,
74 (2011), pp.\ 441--461.

\bibitem{bs} J. F. BONNANS and A. SHAPIRO, {\em Perturbation Analysis of Optimization Problems}, Series in Operations Research, Springer, New York,
2000.

\bibitem{bz} J. M. BORWEIN and Q. J. ZHU, {\em Techniques of Variational Analysis},
CMS Books in Mathematics, Vol.\ 20, Springer, New York, 2005.

\bibitem{ccyy} N. H. CHIEU, T. D. CHUONG, J.-C. YAO and N. D. YEN, {\em Characterizing convexity of a function by its Fr\'echet and limiting second-order
subdifferentials}, Set-Valued Var. Anal., 19 (2011), pp.\ 75--96.

\bibitem{ch} N. H. CHIEU and N. D. HUY, {\em Second-order subdifferentials and convexity of real-valued functions}, Nonlinear Anal., 74 (2011), pp.\ 154--160.

\bibitem{ct}  N. H. CHIEU and N. T. Q. TRANG, {\em Coderivatives and monotonicity of continuouis mappings}, Taiwan. J. Math., to appear (2011).

\bibitem{c} F. H. CLARKE, {\em Optimization and Nonsmooth Analysis}, CMS Books in Mathematics, Vol.\ 1, Wiley, New York, 1983.

\bibitem{chhm} G. COLOMBO, R. HENRION, N. D. HOANG and B. S. MORDUKHOVICH, {\em Optimal control of the sweeping process}, Dynam. Cont. Disc. Impul. Syst., Ser. B, 19 (2012).

\bibitem{dr} A. L. DONTCHEV and R. T. ROCKAFELLAR, {\em Characterizations of strong regularity for variational
inequalities over polyhedral convex sets}, SIAM J. Optim., 6 (1996), pp.\ 1087--1105.

\bibitem{dr-b} A. L. DONTCHEV and R. T. ROCKAFELLAR, {\em Implicit Functions and Solution Mappings:
A View from Variational Analysis}, Springer Monographs in Mathematics, Springer,
Dordrecht, 2009.

\bibitem{ew} A. C. EBERHARD and R. WENCZEL, {\em A study of tilt-stable optimality and sufficient conditions}, Nonlinear Anal., 75 (2012), No. 3.

\bibitem{hmn} R. HENRION, B. S. MORDUKHOVICH and N. M. NAM, {\em Second-order analysis of polyhedral
systems in finite and infinite dimensions with applications to robust stability of variational inequalities},
SIAM J. Optim., 20 (2010), pp.\ 2199--2227.

\bibitem{hos} R. HENRION, J. V. OUTRATA and T. SUROWIEC, {\em On the coderivative of normal cone mappings to inequality systems}, Nonlinear Anal., 71 (2009),
pp.\ 1213--1226.

\bibitem{hos1} R. HENRION, J. V. OUTRATA and T. SUROWIEC, {\em Analysis of M-stationarity points to an EPEC modeling ologopolistic competion in an electricity
spot market}, ESAIM: Contr. Optim. Calc. Var., to appear (2012).

\bibitem{hr} R. HENRION and W. R\"{O}MISCH, {\em On M-stationary points for a stochastic equilibrium problem under equilibrium constraints in
electricity spot market modeling}, Appl. Math., 52 (2007), pp.\ 473--494.

\bibitem{hy} N. Q. HUY and J.-C. YAO, {\em Exact formulae for coderivatives of normal cone mappings to perturbed polyhedral convex sets}, J. Optim. Theory Appl., to appear (2012).

\bibitem{jl} V. JEYAKUMAR and D. T. LUC, {\em Nonsmooth Vector Functions and Continuous Optimization}, Optimization and Its Applications, Vol. 10,
Springer, New York, 2008.

\bibitem{lm} A. B. LEVY and B. S. MORDUKHOVICH, {\em Coderivatives in
parametric optimization}, Math. Program., 99 (2004), pp.\ 311--327.

\bibitem{lpr} A. B. LEVY, R. A. POLIQIUN and R. T. ROCKAFELLAR,
{\em Stability of locally optimal solutions}, SIAM J. Optim., 10 (2000), pp.\  580--604.

\bibitem{lr} A. B. LEVY and R. T. ROCKAFELLAR, {\em Variational conditions and the proto-differentiation of partial subgradient mappings}, Nonlinear Anal., 26 (1996), pp.\
1951--1964.

\bibitem{lz} A. S. LEWIS and S. ZHANG, {\em Partial smoothness, tilt stabiity, and generalized Hessians}, preprint (2011).

\bibitem{m76} B. S. MORDUKHOVICH, {\em Maximum principle
in problems of time optimal control with nonsmooth constraints},
J. Appl. Math. Mech.,  40 (1976), pp.\ 960--969.

\bibitem{m80} B. S. MORDUKHOVICH,  {\em Metric approximations
and necessary optimality conditions for general classes of
extremal problems}, Soviet Math. Dokl., 22 (1980), pp.\ 526--530.

\bibitem{m92} B. S. MORDUKHOVICH,  {\em Sensitivity analysis
in nonsmooth optimization}, in Theoretical Aspects of Industrial
Design, D. A. Field and V. Komkov (eds.), Proceedings in Applied
Mathematics, Vol.\ 58, SIAM, Philadelphia, 1992, pp.\ 32--46.

\bibitem{m93} B. S. MORDUKHOVICH, {\em Complete characterizations of
openness, metric regularity, and Lipschitzian properties of
multifunctions}, Trans. Amer. Math. Soc., 340 (1993), pp.\ 1--35.

\bibitem{m94a} B. S. MORDUKHOVICH, {\em Generalized differential calculus for
nonsmooth and set-valued mappings}, J. Math. Anal. Appl., 183
(1994), pp.\ 250--288.

\bibitem{m94b} B. S. MORDUKHOVICH, {\em Stability theory
for parameteric generalized equations and variational inequalities
via nonsmooth analysis}, Trans. Amer. Math. Soc., 343 (1994), pp.\
609--658.

\bibitem{m02} B. S. MORDUKHOVICH, {\em Calculus of second-order subdifferentials in infinite dimensions}, Control
Cybernet., 31 (2002), pp.\ 557--573.

\bibitem{m06a} B. S. MORDUKHOVICH, {\em Variational Analysis and
Generalized Differentiation, I: Basic Theory}, Grundlehren Series
(Fundamental Principles of Mathematical Sciences), Vol.\ 330,
Springer, Berlin, 2006.

\bibitem{m06b} B. S. MORDUKHOVICH, {\em Variational Analysis and
Generalized Differentiation, II:  Applications}, Grundlehren
Series (Fundamental Principles of Mathematical Sciences), Vol.\
331, Springer, Berlin, 2006.

\bibitem{mo} B. S. MORDUKHOVICH and J. V. OUTRATA, {\em Second-order
subdifferentials and their applications}, SIAM J. Optim., 12
(2001), pp.\ 139--169.

\bibitem{mo1} B. S. MORDUKHOVICH and J. V. OUTRATA, {\em Coderivative analysis
of quasi-variational inequalities with applications to stability in optimization},
SIAM J. Optim., 18 (2007), pp.\ 389--412.

\bibitem{ms} B. S. MORDUKHOVICH and Y. SHAO, {\em Nonsmooth sequential anaysis in Asplund spaces}, Trans. Amer. Math. Soc., 348 (1996), 1235--1280.

\bibitem{mw} B. S. MORDUKHOVICH and B. WANG, {\em Restrictive metric regularity and generalizerd differential calculus in Banach spaces},
Int. J. Maths. Math. Sci., 50 (2004), pp.\ 2650--2683.

\bibitem{nam} N. M. NAM, {\em Coderivatives of the normal cone mappings and Lipschitzian
stability of parametric variational inequalities}, Nonlinear Anal.,
73 (2010), pp.\ 2271--2282.

\bibitem{o} J. V. OUTRATA,  {\em Mathematical programs
with equilibrium constraints: Theory and numerical methods}, in
Nonsmooth Mechanics of Solids, CISM Lecture Notes, Vol.\ 485, J.
Haslinger and G. E. Stavroulakis (eds.), Springer, New York, 2006,
pp.\ 221--274.

\bibitem{pr92} R. A. POLIQUIN and R. T. ROCKAFELLAR, {\em Amenable
functions in optimization}, in Nonsmooth Optimization Methods and
Applications, F. Giannessi (ed.),  Gordon and Breach,
Philadelphia, 1992, pp.\ 338--353.

\bibitem{pr96} R. A. POLIQUIN and R. T. ROCKAFELLAR, {\em Prox-regular
functions in variational analysis}, Trans. Amer. Math. Soc., 348
(1996), pp.\ 1805--1838.

\bibitem{pr98} R. A. POLIQUIN and R. T. ROCKAFELLAR, {\em Tilt stability
of a local minimum}, SIAM J. Optim., 8 (1998), pp.\ 287--299.

\bibitem{q} N. T. QUI, {\em Nonlinear perturbations of polyhedral normal cone
mappings and affine variational inequalities}, J. Optim. Theory Appl., to appear (2012).

\bibitem{rob} S. M. ROBINSON, {\em Strongly regular generalized equations},
Math. Oper. Res., 5 (1980), pp.\ 43--62.

\bibitem{rob1} S. M. ROBINSON, {\em Some continuity properties of polyhedral multifunctions}, Math. Programming Stud., 14 (1981), pp.\ 206--214.

\bibitem{r85} R. T. ROCKAFELLAR, {\em Maximal monotone relations and the second
derivatives of nonsmooth functions}, Ann. Inst. H. Poincar\'e:
Analyse Non Lin\'eaire,  2 (1985), 167--184.

\bibitem{r00} R. T. ROCKAFELLAR, {\em Extended nonlinear
programming}, in Nonlinear Optimization and Related Topics, G. Di
Pillo and F. Giannessi (eds.), Applied Optimization, Vol.\ 36,
Kluwer Academic Publishers, Dordrecht, 2000, pp.\ 381--399.

\bibitem{rw} R. T. ROCKAFELLAR and R. J-B WETS, {\em Variational
Analysis}, Grundlehren Series (Fundamental Principles of
Mathematical Sciences), Vol.\ 317, Springer, Berlin, 2006.

\bibitem{s} W. SCHIROTZEK, {\em Nonsmooth Analysis}, Universitex, Springer, Berlin, 2007.

\bibitem{su} T. SUROWIEC, {\em Explicit Stationarity Conditions and Solution Characterization for Equilibrium Problems with Equilibrium Constraints}, Ph. D. Disseration, Humboldt
University, Berlin, 2010.

\bibitem{yy} J.-C. YAO and N. D. YEN, {\em Coderivative calculation related to a parametric affine variational inequality.
Part~1: Basic calculation}, Acta Math. Vietnam., 34 (2009), pp.\ 157--172.

\bibitem{yy1} J.-C. YAO and N. D. YEN, {\em Coderivative calculation related to a parametric affine variational inequality.
Part~2: Applications}, Pac. J. Optim., 5 (2009), pp.\ 493-?506.

\bibitem{ye} J. J. YE, {\em Constraint qualifications
and necessary optimality conditions for optimization problems with
variational inequality constraints}, SIAM J. Optim., 10 (2000),
pp.\ 943--962.

\bibitem{z} R. ZHANG, {\em Problems of hierarchical optimization in finite dimensions},
SIAM J. Optim., 4 (1994), pp.\ 521--536.
\end{thebibliography}
\end{document}